\documentclass[11pt]{article}
\usepackage{latexsym,amssymb,enumerate,amsmath,epsfig,amsthm,authblk,dsfont}
\usepackage[margin=1in]{geometry}
\usepackage{setspace,color}
\usepackage{tikz}
\usepackage{floatrow}
\usepackage{graphicx,subfigure}
\usepackage[ruled]{algorithm2e}
\usepackage{epstopdf}
\usepackage{grffile}
\usepackage{hyperref}

\newcommand{\p}{\mathbf{p}}
\newcommand{\xb}{\mathbf{x}}
\newcommand{\n}{\mathbf{n}}
\newcommand{\s}{\mathbf{s}}
\newcommand{\D}{\mathbf{D}}
\newcommand{\cof}{\mathrm{cof}}
\newtheorem{thm}{Theorem}[section]
\newtheorem{remark}{Remark}[section]

\begin{document}
\title{Operator Splitting/Finite Element Methods for the Minkowski Problem}
\date{}
\author{
Hao Liu\thanks{Department of Mathematics, Hong Kong Baptist University, Kowloon Tong, Hong Kong. Email: {\bf haoliu@hkbu.edu.hk}.}
,
Shingyu Leung\thanks{Department of Mathematics, The Hong Kong University of Science and Technology, Clear Water Bay, Hong Kong. Email: {\bf masyleung@ust.hk}.}
,
Jianliang Qian\thanks{Department of Mathematics and Department of CMSE, Michigan State University, East Lansing, MI 48824, USA. Email: {\bf jqian@msu.edu}. }
}
\maketitle
\begin{abstract}
	The classical Minkowski problem for convex bodies has deeply influenced the development of differential geometry. During the past several decades, abundant mathematical theories have been developed for studying the solutions of the Minkowski problem, however, the numerical solution of this problem has been largely left behind, with only few methods available to achieve that goal. In this article, focusing on the two--dimensional  Minkowski problem with Dirichlet boundary conditions, we introduce two solution methods, both based on operator--splitting. One of these two methods deals directly with the Dirichlet condition, while the other method uses an approximation \`{a} la Robin of this Dirichlet condition. This relaxation of the Dirichlet condition makes this second method better suited than the first one to treat those situations where the Minkowski equation (of the Monge--Amp\`{e}re type) and the Dirichlet condition are not compatible. Both methods are generalizations of the solution method for the canonical Monge--Amp\`{e}re equation discussed by Glowinski et al. (A Finite Element/Operator--Splitting Method for
	the Numerical Solution of the Two Dimensional Elliptic Monge--Amp\`{e}re Equation, {\it Journal of Scientific Computing}, 79(1), 1--47, 2019); as such they take advantage of a divergence formulation of the Minkowski problem, well--suited to a mixed finite element approximation, and to the the time--discretization via an operator--splitting scheme, of an associated initial value problem. Our methodology can be easily implemented on convex domains of rather general shape (with curved boundaries, possibly). The numerical experiments we performed validate both methods and show that if one uses continuous piecewise affine finite element approximations of the solution of the Minkowski problem and of its three second order derivatives, these two methods provide nearly second order accuracy for the $L^2$ and $L^{\infty}$ norms of the approximation error. (This assumes, of course, that the Minkowski--Dirichlet problem has a smooth solution). One can extend easily the methods discussed in this article, to address the solution of three--dimensional Minkowski problem.

\end{abstract}

\section{Introduction}
\label{sec.intro}
The Minkowski problem (named after Hermann Minkowski (1864--1909)) is an important problem in Differential Geometry. It asks for the construction of a compact surface $S$ as boundary of a convex bounded domain, knowing its Gaussian curvature. 
Given a compact strictly convex hypersurface $S$ in the $d$--dimensional real space $\mathds{R}^{d}$, the Gauss map $\mathbf{G}$ is a diffeomorphism from $S$ to the unit sphere $\mathbf{S}^{d-1}$ of $\mathds{R}^{d}$. Map $\mathbf{G}$ is defined by $\mathbf{G}(\xb)=\n(\xb),\forall \xb\in S$, where $\n(\xb)$ denotes the unit outward normal of $S$ at $\xb$. Accordingly, the Gauss-Kronecker curvature $K$ is the Jacobian of the Gauss map. Minkowski stated that one has
\begin{equation}
	\int_{\mathbf{S}^{n-1}} \xb (K(\mathbf{G}^{-1}(\xb)))^{-1}d\sigma(\xb)=\mathbf{0},
	\label{eq.mink.cond}
\end{equation}
where $\sigma$ is the Lebesgue measure on $\mathbf{S}^{d-1}$. Conversely, Minkowski posed the following (inverse) problem: Suppose that $f$ is a strictly positive function defined over $\mathbf{S}^{d-1}$ verifying $\int_{\mathbf{S}^{n-1}} \xb f(\xb)d\sigma(\xb)=\mathbf{0}$; can one find a hypersurface $S$ having $1/f$ as Gaussian curvature? 
In \cite{minkowski1989allgemeine,minkowski1989volumen}, Minkowski discussed the existence and uniqueness of solutions to the above inverse problem. For $d=2$, the solution regularity was proved by Lewy \cite{lewy1935priori,lewy1938differential}, Nirenberg \cite{nirenberg1953weyl}, and Pogorelov \cite{pogorelov1952regularity}, while, for $d>2$,
the solution regularity was analyzed by Cheng and Yau \cite{cheng1976regularity} and Pogorelov \cite{pogorelov1978minkowski}.

Despite being around for more than a century and being one of the most important problems in Differential Geometry, no much was done concerning the numerical solution of the Minkowski problem. The earliest attempt we could find was discussed in \cite{little1983iterative,little1985recovering}, two publications dedicated to the solution of a related problem: namely, reconstructing a shape from extended Gaussian images. In \cite{lamberg1995minkowski}, after generalizing Minkowski's proof, Lamberg converted the Minkowski problem into an optimization one, the resulting algorithm solving a polyhedral version of the Minkowski problem. In \cite{lamberg2001numerical}, Lamberg introduced an algorithm based on Minkowski's isoperimetric inequality, leading to an approximate Minkowski problem taking place in a finite-dimensional function space spanned by truncated spherical harmonic series. In a more recent publication \cite{cheng2003construction}, Cheng designed a level-set based finite-difference PDE method to drive an implicitly defined surface towards shapes arising from the Minkowski problem.

In all the above cited works the hypersurface is supposed to be closed. Actually, another type of Minkowski problem is the Minkowski--Dirichlet problem. For the Minkowski--Dirichlet problem, one supposes that the hypersurface is open and bounded, with a Dirichlet condition imposed on its boundary. The well--posedness of this problem has been addressed by many authors: For example, Bakelman \cite{bakel1980dirichlet}, Lions \cite{lions1985equations} and Urbas \cite{urbas1984equation,urbas1988global,urbas1991boundary} have proved the existence and uniqueness of a solution. Trudinger and Urbas \cite{trudinger1983dirichlet} proved a necessary and sufficient condition for the classical solvability of the Minkowski--Dirichlet problem. Recently, in \cite{hamfeldt2018convergent} Hamfelt designed  a monotone finite-difference method to solve the Minkowski--Dirichlet problem; since the method relies on wide stencils, it is advantageous for those situations where, due to the lack of classical solutions, one looks for viscosity solutions.

Here, we propose two new methods for the numerical solution of the Minkowski--Dirichlet problem. The first method, well suited to problems with classical solutions, imposes the Dirichlet condition in a strong sense. On the other hand, the second method imposes the Dirichlet condition in a least--squares sense (via a quadratic penalty technique), making it appropriate for those situations where, due to data incompatibility, the Minkowski--Dirichlet problem has no solution. Of course, the second method has also the ability to capture classical solutions, if such solutions do exist. The Minkowski problem we will look at can be described as follows: Let $\Omega$ be a bounded domain of $\mathds{R}^d$ and $K$ be a positive function defined over $\Omega$, and let $g$ be a function defined on the boundary $\partial\Omega$; can one find a function $u$ defined over $\Omega$ and verifying $u|_{\partial \Omega}=g$, such that $K$ is the Gauss curvature of the graph of $u$ (a surface in $\mathds{R}^{d+1}$)? In partial differential equation form, the above Minkowski--Dirichlet problem reads as follows:
\begin{equation}
	\begin{cases}
		\frac{\det(\D^2u)}{(1+|\nabla u|^2)^{1+d/2}}=K &\mbox{ in } \Omega,\\
		u=g &\mbox{ on } \partial \Omega.
	\end{cases}
	\label{eq.mk}
\end{equation}
The partial differential equation in (\ref{eq.mk}) belongs to a family of Monge--Amp\`{e}re equations. The simplest element of this family is clearly the following canonical Monge--Amp\`{e}re equation 
\begin{equation}
	\det(\D^2 u)=f \mbox{ in } \Omega.
	\label{eq.MA}
\end{equation}
Equation (\ref{eq.MA}) is elliptic if $f>0$. The above Monge--Amper\`{e}re equation (\ref{eq.MA}) is a fully nonlinear second order partial differential equation; it has been drawing a lot of attention lately, mostly because its relations with optimal transport problems (other applications are described in, e.g., \cite{feng2013recent}; see also the references therein). During the past three decades, a variety of methods have been designed to solve numerically equation (\ref{eq.MA}), completed by boundary conditions (mostly Dirichlet's) (some of these methods are described in the review article \cite{feng2013recent}). As expected, most of these methods focus on the two--dimensional Monge--Amp\`{e}re equation and cover a large variety of approaches. Combinations of (mixed) finite element approximations and augmented Lagrangian or least-squares formulations have been applied to the solution of (\ref{eq.MA}) and related fully nonlinear elliptic equations such as Pucci's (see \cite{benamou2000computational,caboussat2013least,dean2003numerical,glowinski2006augmented,dean2004numerical,dean2008numerical,dean2005operator,dean2006numerical,glowinski2015variational,glowinski2017splitting,glowinski2009numerical,mohammadi2007optimal,caffarelli2008numerical,feng2013recent} for details). Alternative finite-difference and finite-element methods have been developed for these fully nonlinear elliptic equations as well; see  \cite{awanou2014standard,benamou2010two,brenner2011c,brenner2012finite,froese2011fast,froese2011convergent,loerap05,schaeffer2016accelerated,glowinskillq2017,glowinskillq20173d,feng2013recent}, this list being far from complete.

The main goal of this article is to extend to problem (\ref{eq.mk}) (assuming $d$ = 2), the operator--splitting based methods developed in \cite{glowinskillq2017,glowinskillq20173d} for the solution of equation (\ref{eq.MA}) (completed by Dirichlet conditions) in dimensions 2 and 3 and in \cite{glowinski2020numerical,liu2022efficient} for the eigenvalue problems of (\ref{eq.MA}). Following \cite{glowinskillq2017,glowinskillq20173d}, the first step in that direction is to take advantage of a divergence formulation of problem (\ref{eq.mk}), better suited to finite element approximations. The second step is to decouple (in some sense) differential operators and nonlinearities by introducing as additional unknown functions $\p=\D^2u$ (as done in \cite{glowinski2020numerical,glowinskillq20173d}) and $\s=\nabla u$ (which was not necessary in \cite{glowinski2020numerical,glowinskillq20173d}). At the end of the second step, one has replaced the highly nonlinear scalar Minkowski equation by an equivalent system of linear and nonlinear equations for $u,\p$ and $\s$, whose formalism is simpler.  In the third step, we associate an initial value problem (IVP) with the above system and use operator--splitting to time--discretize the above IVP, in order to capture its steady state solution(s). We use simple finite-element approximations of the mixed type to implement the above methodology: indeed, we use finite-element spaces of continuous piecewise affine functions to approximate $u$ and its three second-order derivatives, making our methods well--suited to solve problem (\ref{eq.mk}) on domains $\Omega$ with curved boundaries.

As mentioned above we will develop two new methods for the solution of problem (\ref{eq.mk}): these two methods are very close to each other, the first one dealing directly with the boundary condition $u = g$ on $\partial \Omega$, while the second one imposing the boundary condition in a least--squares sense.   

This article is organized as follows: In Section \ref{sec.problem}, we state some theoretical results on the existence and uniqueness of solutions to the Minkowski--Dirichlet problem (\ref{eq.mk}). In Section \ref{sec.div}, we provide the divergence formulation of problem (\ref{eq.mk}) and associate with it two initial value problems, which differ by the way the Dirichlet boundary condition is treated. The time discretization of these two initial value problems by operator--splitting is discussed in Section \ref{sec.time}, followed by their finite-element space discretization in Section \ref{sec.space}. We address in Section \ref{sec.initial} the initialization of the two above algorithms. In Section \ref{sec.example}, we report the results of numerical experiments validating our methodology. Section \ref{sec.conclusion} concludes the article.

\section{Problem formulation, existence, uniqueness and regularity results}
\label{sec.problem}
We defined the Minkowski problem in Section \ref{sec.intro}. In this article, we will focus on the numerical solution of the Minkowski--Dirichlet problem (\ref{eq.mk}), assuming that $d = 2$ (2-D). A first step to that goal is to rewrite (\ref{eq.mk}) as 
\begin{equation}
	\begin{cases}
		\det(\D^2u)=K(1+|\nabla u|^2)^{1+d/2} & \mbox{ in }\Omega,\\
		u=g &\mbox{ on } \partial\Omega,
	\end{cases}
	\label{eq.mkn}
\end{equation}
a Monge--Amp\`{e}re type formulation, better suited for numerical solution. In (\ref{eq.mkn}), $K$ $(>0)$ is the prescribed curvature and $\D^2u=\left(\frac{\partial^2 u}{\partial x_i\partial x_j}\right)_{1\leq i,j \leq d}$ is the Hessian matrix of function $u$.

To put our computational investigations into perspective, we recall some classical results concerning the existence, uniqueness and regularity of classical solutions to problem (\ref{eq.mkn}). In \cite{trudinger1983dirichlet}, one proves the following results about existence and uniqueness.

\begin{thm}
	Suppose that in (\ref{eq.mkn}), $\Omega$ is a uniformly bounded convex domain of $\mathds{R}^d$, its boundary $\partial\Omega$ having $C^{1,1}$--regularity. Then, problem (\ref{eq.mkn}) has, for any $g\in C^{1,1}(\bar{\Omega})$, a unique solution in $C^2(\Omega)\cap C^{0,1}(\bar{\Omega})$, if and only if 	
	\begin{equation}
		\int_{\Omega} Kdx<\omega_d,
		\label{eq.exist.1}
	\end{equation}
	and
	\begin{equation}
		K=0 \mbox{ on }\partial\Omega.
		\label{eq.exist.2}
	\end{equation}
	The constant $\omega_d$ in (\ref{eq.exist.1}) is given by $\omega_d=\int_{\mathds{R}^d} \frac{d\xi}{(1+|\xi|^2)^{1+d/2}}$ (implying $\omega_2=\pi$ and $\omega_3=4\pi/3$); actually, $\omega_d$ is the volumn of the unit ball of $\mathds{R}^d$.
\end{thm}

Condition (\ref{eq.exist.2}) is required to make sure that a solution exists for arbitrary $g$. It is proved in \cite{trudinger1983dirichlet} that if $K$ does not vanish on the boundary, one can find a smooth function $g$ such that problem (\ref{eq.mkn}) has no solution. This situation is well illustrated by the elementary example in Remark \ref{rmk.counterex}. 
\begin{remark} \label{rmk.counterex}
	We give an example in which problem (\ref{eq.mkn}) has no classical solutions in general if condition (\ref{eq.exist.2}) is not verified.  The particular case we consider is 
	\begin{align}
		\det \D^2u=K(1+|\nabla u|^2)^2 \mbox{ in } \Omega,
		\label{eq.exist.ex}
	\end{align}
	where $K$ is  a constant verifying $0<K<1$, and $\Omega=\left\{(x_1,x_2)\in \mathds{R}^2,x_1^2+x_2^2<1\right\}$. Condition (\ref{eq.exist.1}) is clearly verified. All solutions of equation (\ref{eq.exist.ex}) are either convex or concave and given by
	\begin{align}
		u(x_1,x_2)=c\pm \sqrt{R^2-(x_1-a)^2-(x_2-b)^2}, \forall (x_1,x_2)\in \bar{\Omega},
		\label{eq.exist.ex.1}
	\end{align}
	with $R=1/\sqrt{K}$ $(>1)$, $(a,b,c)\in \mathds{R}^3$, and $R-1\geq \sqrt{a^2+b^2}$. It follows from (\ref{eq.exist.ex.1}) that the traces on $\partial\Omega$ of the solutions to equation (\ref{eq.exist.ex}) belong to the set
	\begin{align*}
		G=\Big\{ g\, |\, &g(x_1,x_2)=c\pm \sqrt{R^2-(x_1-a)^2-(x_2-b)^2}, \; \forall (x_1,x_2)\in \partial\Omega,\\ 
		&(a,b,c)\in \mathds{R}^3,\; \ R-1\geq \sqrt{a^2+b^2}\Big\}.
	\end{align*}
	Suppose now that $g$ is a smooth function not belonging to $G$ (there are infinitely many such functions); it cannot be the trace of a solution to equation (\ref{eq.exist.ex}), showing that in general the related problem (\ref{eq.mkn}) has no solution.                           
\end{remark}

In \cite{urbas1984equation,urbas1988global,urbas1991boundary}, one discusses regularity of the solution in the critical case defined by
\begin{equation}
	\int_{\Omega} Kdx=\omega_d,
	\label{eq.exist.limit}
\end{equation}
where the following results are proved.

\begin{thm}
	Let $\Omega$ be a uniformly convex domain of $\mathds{R}^d$ with a $C^{2,1}$ smooth boundary, and $K$ be a positive $C^2$ smooth function verifying (\ref{eq.exist.limit}). If $u$ is a solution of the Minkowski--Dirichlet (\ref{eq.mkn}), then
	\begin{enumerate}[(i)]
		\item $u\in C^{0,1/2}(\Omega)$;
		\item the graph of $u$ is $C^{2,\alpha}$--smooth for some $\alpha\in(0,1)$;
		\item $u|_{\partial\Omega}$ is $C^{1,\alpha}$--smooth;
		\item if $\partial\Omega$ is $C^{k+1,\alpha}$ and $K\in C^{k-1,\alpha}$ with $k\geq2$, then the graph of $u$ is $C^{k+1,\alpha}$--smooth and $u|_{\partial\Omega}$ is $C^{k+1,\alpha}$--smooth.
	\end{enumerate}
\end{thm}
See \cite{trudinger2008monge} for more details on the solution of the Minkowski problem.

Some of the conditions in the above two theorems are rather restrictive and/or not easy to verify. Nevertheless, the results they are reporting are very useful from two perspectives: on one hand, they suggest test problems, where we know in advance that solutions exist; on the other hand, they also suggest some other examples, where the answer to existence will be indicated by the results of our computations. Finally, we will also consider test problems with known solutions so as to check how fast and how accurately our methods recover them.

\section{Divergence formulations of the 2-D Minkowski problem and relaxation by penalty of the Dirichlet condition}
\label{sec.div}

\subsection{Synopsis}
There are cases where the data $K$ and $g$ do not allow the existence of classical smooth solutions to problem (\ref{eq.mkn}). In \cite{hamfeldt2018convergent}, 
one introduces a notion of viscosity solution to problem (\ref{eq.mkn}), with the solution satisfying the generalized Monge-Amp\`{e}re equation in \cite{bakelman1986generalized}, but not necessarily the Dirichlet condition. In the following sub-sections, we will consider two divergence formulations of problem (\ref{eq.mkn}) in dimension two to enforce the Dirichlet condition. The first formulation keeps the Dirichlet condition as it is and is well-suited to those situations where problem (\ref{eq.mkn}) has classical solutions. On the other hand, the second formulation makes use of penalty to relax the Dirichlet condition; for large values of the penalty parameter, one recovers accurately classical solutions if such solutions do exist, or generalized solutions in the absence of classical solutions.

\subsection{A first divergence formulation of the 2-D Minkowski-Dirichlet problem}
\label{sec.scheme1}
If $d=2$, problem (\ref{eq.mkn}) enjoys the following equivalent formulation (in the sense of distributions): 
\begin{equation}
	\begin{cases}
		-\nabla\cdot\left(\cof(\D^2u)\nabla u\right)+2K(1+|\nabla u|^2)^2=0 &\mbox{ in }\Omega,\\
		u=g & \mbox{ on } \partial\Omega,
	\end{cases}
	\label{eq.mkdiv}
\end{equation}
where matrix $\cof(\D^2u)$ is the cofactor matrix of Hessian $\D^2u$, that is
$$\cof(\D^2u)=
\begin{pmatrix}
	\frac{\partial^2 u}{\partial x_2^2} & -\frac{\partial^2 u}{\partial x_1 \partial x_2}\\
	-\frac{\partial^2 u}{\partial x_1 \partial x_2} & \frac{\partial^2 u}{\partial x_1^2}
\end{pmatrix}.
$$
Problem (\ref{eq.mkdiv}) is equivalent to
\begin{equation}
	\begin{cases}
		\begin{cases}
			-\nabla\cdot\left(\cof(\p)\nabla u\right)+2K(1+|\s|^2)^2=0 &\mbox{ in }\Omega,\\
			u=g & \mbox{ on } \partial\Omega,
		\end{cases}\\
		\p-\D^2u=\mathbf{0} \quad \mbox{ in } \Omega,\\
		\s=\nabla u \quad \mbox{ in } \Omega.
	\end{cases}
	\label{eq.mkdiv.1}
\end{equation}
In order to avoid possible troubles at those points of $\bar{\Omega}$ where $K$ may vanish, we approximate
system (\ref{eq.mkdiv.1}) by
\begin{equation}
	\begin{cases}
		\begin{cases}
			-\nabla\cdot\left(\left(\varepsilon\mathbf{I}+\cof(\p)\right)\nabla u\right)+2K(1+|\s|^2)^2=0 & \mbox{ in } \Omega,\\
			u|_{\partial\Omega}=g & \mbox{ on } \Omega,
		\end{cases}\\
		\p-\D^2u=\mathbf{0},\\
		\s-\nabla u=\mathbf{0},
	\end{cases}
	\label{eq.split}
\end{equation}
with $\varepsilon$ a small positive parameter (in practice we will take $\varepsilon$ of the order of $h^2$, $h$ being a space discretization step). We used successfully this type of regularization in \cite{glowinskillq2017}, for the solution of the canonical Monge-Amp\`{e}re equation (\ref{eq.MA}) completed by Dirichlet boundary conditions. 

To solve system (\ref{eq.split}) we are going to associate with it the following initial value problem 
\begin{equation}
	\begin{cases}
		\begin{cases}
			\frac{\partial u}{\partial t}-\nabla\cdot\left(\left(\varepsilon\mathbf{I}+\cof(\p)\right)\nabla u\right)+2K(1+|\s|^2)^2=0 &\mbox{ in } \Omega\times (0,+\infty),\\
			u|_{\partial\Omega}=g & \mbox{ on } \partial\Omega\times (0,+\infty),
		\end{cases}\\
		\frac{\partial \p}{\partial t}+\gamma_1\left(\p-\D^2u\right)=\mathbf{0} \mbox{ in } \Omega\times (0,+\infty),\\
		\frac{\partial \s}{\partial t}+\gamma_2\left(\s-\nabla u\right)=\mathbf{0} \mbox{ in } \Omega\times (0,+\infty),\\
		(u(0),\p(0),\s(0))=(u_0,\p_0,\s_0),
	\end{cases}
	\label{eq.split1}
\end{equation}
to be time-discretized by operator-splitting (in Section \ref{sec.timedis}). In (\ref{eq.split1}), $\gamma_1$ and $\gamma_2$ are two positive coefficients chosen so that the smooth modes of $\p$ and $\s$ evolve in time roughly at the same speed as that of $u$.  Following \cite{glowinskillq2017}, we advocate defining $\gamma_1$ and $\gamma_2$ by
\begin{eqnarray*}
	\gamma_1&=&\beta_1\lambda_0\left(\varepsilon+\sqrt{\alpha}\right),\\
	\gamma_2&=&\beta_2\lambda_0\left(\varepsilon+\sqrt{\alpha}\right),
\end{eqnarray*}
where $\lambda_0$ is the smallest eigenvalue of operator $-\nabla^2$ in $H_0^1(\Omega)$, $\alpha$ is the lower bound of $K$, and $\beta_1$ and $\beta_2$ are two constants of order one.

We comment in passing that we have used and will continue to use the notation $\phi(t)$ for the function $x\rightarrow \phi(x,t)$.
In Section \ref{sec.initial}, we will discuss the initialization of system (\ref{eq.split1}).

\subsection{A divergence formulation of the 2-D Minkowski-Dirichlet problem with relaxation of the boundary condition}
\label{sec.scheme2}
Remark \ref{rmk.counterex} in Section \ref{sec.problem} gives an example to demonstrate that problem (\ref{eq.mkn}) may have no solution, unless function $g$ belongs to a very specific class of functions. In order to deal with such no-solution scenarios as well as we could, we are going to relax the boundary condition $u = g$ using a penalty technique of the least-squares type. If problem (\ref{eq.mkn}) has a classical solution, we expect to recover it when the penalty parameter converges to $+\infty$.

The simplest way to proceed is to start from the following variational formulation verified (formally) by any solution u of problem (\ref{eq.mkn}): 
\begin{equation}
	\begin{cases}
		u\in H^1(\Omega),\\
		\displaystyle\int_{\Omega} (\cof (\D^2u) \nabla u) \cdot \nabla v dx+ 2\int_{\Omega}K(1+|\nabla u|^2)^2\,v\,dx=0, \ \forall v\in H^1_0(\Omega),\\
		u=g \mbox{ on } \partial\Omega.
	\end{cases}\\
	\label{eq.mkma.var}
\end{equation}
In order to relax the Dirichlet boundary condition, we are going to apply to problem (\ref{eq.mkma.var}) the well-known penalty method discussed in  \cite{glowinski1984numerical,glowinski2008lectures} to approximate Dirichlet's problems for linear second-order elliptic operators by Robin's ones.  

Let $\kappa$ be a positive constant. We (formally) approximate the variational problem (\ref{eq.mkma.var}) by 
\begin{equation}
	\begin{cases}
		u\in H^1(\Omega),\\
		\displaystyle\int_{\Omega} (\cof(\D^2u)\nabla u) \cdot \nabla v dx +2\int_{\Omega}K(1+|\nabla u|^2)^2vdx\; +\\
		\quad\; \kappa \int_{\partial\Omega} (u-g)vd\Gamma =0,\ \forall v\in H^1(\Omega), 
	\end{cases}
	\label{eq.mkma.varpenalty}
\end{equation}
where coefficient $\kappa$ acts as a weight, controlling the level of penalization. Some remarks are in order. 

\begin{remark}
	Let us consider the functional $j_2: H^1(\Omega)\rightarrow \mathds{R}$ defined by
	\begin{equation*}
		j_2(u)=\frac{\kappa}{2}\int_{\partial \Omega} |u-g|^2 d\Gamma, \ \forall v\in H^1(\Omega).
	\end{equation*}
	Functional $j_2$ is convex and $C^{\infty}$ over $H^1(\Omega)$, its differential $Dj_2(v)$ at $v$ being given by
	\begin{align}
		\langle D j_2(v),w\rangle =\kappa\int_{\partial\Omega} (v-g)w d\Gamma, \ \forall v,w\in H^1(\Omega), 
		\label{eq.j2.diff}
	\end{align}
	where $\langle \cdot,\cdot\rangle$ denotes a duality pairing between $(H^1(\Omega))'$ (the dual space of $H^1(\Omega)$) and $H^1(\Omega)$. Consequently, we can identify $Dj_2(u)$ with $\kappa(u|_{\partial\Omega}-g)$ and replace $\kappa\int_{\partial\Omega} (u-g)vd\Gamma$ in  (\ref{eq.mkma.varpenalty}) by $\langle Dj_2(u),v\rangle$. 
\end{remark}

\begin{remark}
	If a function $u$ is a solution of the nonlinear variational problem (\ref{eq.mkma.varpenalty}), it is also a solution (in the sense of distributions) of the following (fully nonlinear) boundary value problem  
	\begin{align}
		\begin{cases}
			-\nabla\cdot (\cof(\D^2u)\nabla u) +2K(1+|\nabla u|^2)^2=0 & \mbox{ in } \Omega,\\
			\frac{1}{\kappa} (\cof(\D^2u)\nabla u) \cdot \n +u=g & \mbox{ on } \partial \Omega,
		\end{cases}
		\label{eq.var.n}
	\end{align}
	where, in (\ref{eq.var.n}), $\n$ denotes the unit outward normal vector at $\partial\Omega$. The boundary condition in (\ref{eq.var.n}) is a (nonlinear) \emph{Robin} boundary condition. When $\kappa\rightarrow +\infty$, problem (\ref{eq.var.n}) `converges' (formally) to problem (\ref{eq.mkn}), justifying our second divergence formulation of problem (\ref{eq.mkn}).
\end{remark}

\begin{remark} \label{rmk.l1}
	A natural alternative to problem (\ref{eq.mkma.varpenalty}) is the one described by 
	\begin{equation}
		\begin{cases}
			u\in H^1(\Omega),\\
			\displaystyle\int_{\Omega} (\cof(\D^2u)\nabla u) \cdot \nabla v dx +2\int_{\Omega}K(1+|\nabla u|^2)^2vdx\; +\\
			\quad\; \langle \partial j_1(u),v \rangle =0,\ \forall v\in H^1(\Omega),
		\end{cases}
		\label{eq.mkma.varpenalty.l1}
	\end{equation}
	where, in (\ref{eq.mkma.varpenalty.l1}), $\partial j_1(u)$ is the sub-differential at $u$ of the convex Lipschitz continuous functional $j_1:H^1(\Omega)\rightarrow \mathds{R}$, defined by
	\begin{align*}
		j_1(v)=\kappa\int_{\partial\Omega} |v-g|d\Gamma, \ \forall v\in H^1(\Omega).
	\end{align*}
	This type of $L^1$ functional is very common in \emph{Non-Smooth Mechanics} and increasingly popular in \emph{Data Science} as shown by various chapters of \cite{glowinski2017splitting}.
\end{remark}

Proceeding as in Section \ref{sec.scheme1}, we associate with (\ref{eq.mkma.varpenalty}) the semi-variational system 
\begin{equation}
	\begin{cases}
		\begin{cases}
			u\in H^1(\Omega),\\
			\int_{\Omega} \left((\varepsilon \mathbf{I} +\cof(\p))\nabla u\right) \cdot \nabla v dx\, +\, \int_{\Omega} \varepsilon\nabla u \cdot \nabla vdx\;+\\
			\quad\quad\;  +\, 2\int_{\Omega}K(1+|\s|^2)^2vdx \ + \ \kappa \int_{\partial\Omega}(u-g)vd\Gamma=0,\quad \forall v\in H^1(\Omega),
		\end{cases}\\
		\p-\D^2 u=\mathbf{0},\\
		\s-\nabla u=\mathbf{0}.
	\end{cases}
	\label{eq.penalty}
\end{equation}
The next step is to associate with (\ref{eq.penalty}) an initial value problem, as we have done with (\ref{eq.split}) in Section \ref{sec.scheme1}. The initial value problem reads as:
\begin{equation}
	\begin{cases}
		\begin{cases}
			u(t)\in H^1(\Omega), \; \forall t>0,\\
			\displaystyle\int_{\Omega} \frac{\partial u}{\partial t}v dx\;+\;\int_{\Omega} [\varepsilon I+\cof (\p)]\nabla u\cdot \nabla v dx \\
			\quad\;+\; 2\displaystyle\int_{\Omega}K(1+|\s|^2)^2vdx +\kappa\int_{\partial\Omega} (u-g)vd\Gamma=0,\quad \forall v\in H^1(\Omega),
		\end{cases}\\
		\frac{\partial \p}{\partial t} +\gamma_1(\p-\D^2u)=\mathbf{0} \mbox{ in } \Omega\times (0,+\infty),\\
		\frac{\partial \s}{\partial t}+\gamma_2(\s-\nabla u)=\mathbf{0} \mbox{ in } \Omega \times (0,+\infty),\\
		(u(0),\p(0),\s(0))=(u_0,\p_0,\s_0).
	\end{cases}
	\label{eq.split2}
\end{equation}
As in Section \ref{sec.scheme1}, we advocate taking
\begin{eqnarray*}
	\gamma_1&=&\beta_1\lambda_0\left(2\varepsilon+\sqrt{\alpha}\right),\\
	\gamma_2&=&\beta_2\lambda_0\left(2\varepsilon+\sqrt{\alpha}\right).
\end{eqnarray*}
In Section \ref{sec.initial}, we will discuss the initialization of system (\ref{eq.split2}).

\section{Discretization of the IVPs (\ref{eq.split1}) and (\ref{eq.split2}) by operator-splitting}
\label{sec.time}

In this section, we are going to apply the Lie scheme to the time-discretization of the initial value problems (\ref{eq.split1}) and (\ref{eq.split2}); see \cite{glowinski2017splitting} for details on the Lie scheme.

In the following, let $\Delta t$ $(>0)$ denote a time-discretization step, $t^n=n\Delta t$, and let $(u^n,\p^n,\s^n)$ denote an approximation of $(u,\p,\s)$ at $t=t^n$.

\subsection{Time discretization of the initial value problem (\ref{eq.split1})}
\label{sec.timedis}

The Lie-scheme we employ here is a variant of the one we used in \cite{glowinskillq2017} to solve the Monge-Amp\`{e}re equation (\ref{eq.mk}) completed by a Dirichlet boundary condition. It reads as: 
\begin{align}
	(u^0,\p^0,\s^0)=(u_0,\p_0,\s_0).
\end{align}
For $n\geq 0$, $(u^n,\p^n,\s^n)\rightarrow (u^{n+1/2},\p^{n+1/2},\s^{n+1/2}) \rightarrow (u^{n+1},\p^{n+1},\s^{n+1})$ as follows:

\noindent\underline{\emph{The First Fractional Step:}}\\
Solve
\begin{equation}
	\begin{cases}
		\begin{cases}
			\frac{\partial u}{\partial t}-\nabla\cdot\left[\left(\varepsilon\mathbf{I}+\cof(\p^n)\right)\nabla u\right]+
			2K(1+|\s^n|^2)^2=0 & \mbox{ in } \Omega\times(t^n,t^{n+1}) ,\\
			u=g &\mbox{ on } \partial \Omega\times (t^n,t^{n+1}) ,
		\end{cases}\\
		\frac{\partial \p}{\partial t}=\mathbf{0} \mbox{ in } \Omega\times(t^n,t^{n+1})  ,\\
		\frac{\partial \s}{\partial t}=\mathbf{0} \mbox{ in } \Omega\times(t^n,t^{n+1})  ,\\
		(u,\p,\s)(t^n)=(u^n,\p^n,\s^n),
	\end{cases}
	\label{eq.split1.1}
\end{equation}
and set
\begin{equation}
	u^{n+1/2}=u(t^{n+1}),\quad \p^{n+1/2}=\p(t^{n+1})(=\p^n),\quad \s^{n+1/2}=\s(t^{n+1})(=\s^n).
	\label{eq.split1.2}
\end{equation}

\underline{\emph{The Second Fractional Step:}}\\
Solve
\begin{equation}
	\begin{cases}
		\frac{\partial u}{\partial t}=0 &\mbox{ in } \Omega\times(t^n,t^{n+1}) ,\\
		\frac{\partial \p}{\partial t}+\gamma_1(\p-\mathbf{D}^2 u^{n+1/2})=0 &\mbox{ in } \Omega\times(t^n,t^{n+1}),\\
		\frac{\partial \s}{\partial t}+\gamma_2(\s-\nabla u^{n+1/2})=0 & \mbox{ in } \Omega\times(t^n,t^{n+1}),\\
		(u,\p,\s)(t^n)=(u^{n+1/2},\p^{n+1/2},\s^{n+1/2}) ,
	\end{cases}
	\label{eq.split1.3}
\end{equation}
and set
\begin{equation}
	u^{n+1}=u(t^{n+1})(=u^{n+1/2}),\quad \p^{n+1}=P_+\left[\p(t^{n+1})\right],\quad \s^{n+1}=\s(t^{n+1}).
	\label{eq.split1.4}
\end{equation}

In (\ref{eq.split1.4}), $P_+(\cdot)$ is a (kind of) projection operator which maps the space of the $2\times 2$ symmetric matrices onto the closed cone of the $2\times 2$ symmetric positive semi-definite matrices; we will return to operator $P_+$ in Section \ref{sec.P}.

We still need to solve the initial value problems that one encounters in (\ref{eq.split1.1}) and (\ref{eq.split1.3}). There is no difficulty with (\ref{eq.split1.3}) since the three initial value problems it contains have closed form solutions, leading to
$$
\begin{cases}
	u(n+1)=u^{n+1/2},\\
	\p(t^{n+1})=e^{-\gamma_1\Delta t}\p^n+\left(1-e^{-\gamma_1\Delta t}\right) \mathbf{D}^2 u^{n+1},\\
	\s(t^{n+1})=e^{-\gamma_2\Delta t}\s^n+\left(1-e^{-\gamma_2\Delta t}\right) \nabla u^{n+1}.
\end{cases}
$$

It remains to solve the parabolic problem (\ref{eq.split1.1}); for its solution, we advocate performing just one step of the backward Euler scheme. We obtain then
$$
\begin{cases}
	\frac{u^{n+1}-u^n}{\Delta t}-\nabla\cdot\left[\left(\varepsilon\mathbf{I}+
	\cof(\p^n)\right)\nabla u^{n+1}\right]+2K(1+|\s^n|^2)^2=0 &\mbox{ in } \Omega \, ,\\
	u^{n+1}=g &\mbox{ on } \partial\Omega,
\end{cases}
$$
a (relatively) simple Dirichlet problem for a linear self-adjoint second-order strongly elliptic operator with variable coefficients, well-suited to finite-element approximations as we shall see in Section \ref{sec.space}.  

Collecting the above results, we will employ the following time-discretization scheme to solve the initial value problem (\ref{eq.split1}):
\begin{equation}
	(u^0,\p^0,\s^0)=(u_0, \p_0,\s_0).
	\label{eq.split1Dis.0}
\end{equation}
For $n\geq 0$, $(u^n,\p^n,\s^n)\rightarrow (u^{n+1},\p^{n+1},\s^{n+1})$ as follows:

\noindent Solve
\begin{eqnarray}
	&&\begin{cases}
		\frac{u^{n+1}-u^n}{\Delta t}-\nabla\cdot\left[\left(\varepsilon\mathbf{I}+
		\cof(\p^n)\right)\nabla u^{n+1}\right]+2K(1+|\s^n|^2)^2=0 &\mbox{ in } \Omega,\\
		u^{n+1}=g &\mbox{ on } \partial\Omega,
	\end{cases}
	\label{eq.split1Dis.1}
\end{eqnarray}
and compute
\begin{align}
	\begin{cases}
		\p^{n+1}=P_+\left[e^{-\gamma_1\Delta t}\p^n+
		\left(1-e^{-\gamma_1\Delta t}\right) \mathbf{D}^2 u^{n+1}\right], \\
		\s^{n+1}=e^{-\gamma_2\Delta t}\s^n+
		\left(1-e^{-\gamma_2\Delta t}\right) \nabla u^{n+1} . 
	\end{cases}
	\label{eq.split1Dis.2}
\end{align}

\subsection{Time discretization of the initial value problem (\ref{eq.split2})}
As expected, there are many commonalities between the ways we discretize systems (\ref{eq.split1}) and (\ref{eq.split2}); we will take advantage of them. The major difference is that it is much easier to operate directly on the variational formulation of the Monge-Amp\`{e}re part of the problem so as to avoid dealing explicitly with the complicated Robin condition we visualized in (\ref{eq.var.n}). Denote the updated boundary condition at $t^n$ by $g^n$. The Lie scheme we are going to use reads as:  
\begin{eqnarray}
	&&\quad (u^0,\p^0,\s^0,g^0)=(u_0,\p_0,\s_0,g).	\label{eq.split2.0} \\
	\mbox{For }\; n\geq 0,&&\quad (u^n,\p^n,\s^n,g^n)\rightarrow (u^{n+1/3},\p^{n+1/3},\s^{n+1/3},g^{n+1/3})\rightarrow \nonumber\\ 
	&&\rightarrow (u^{n+2/3},\p^{n+2/3},\s^{n+2/3},g^{n+2/3})\rightarrow (u^{n+1},\p^{n+1},\s^{n+1},g^{n+1}),	\nonumber
\end{eqnarray}
where we outline the three fractional steps as the following.

\noindent\underline{\emph{The First Fractional Step:}}\\
Solve
\begin{equation}
	\begin{cases}
		\begin{cases}
			u(t)\in H^1(\Omega), \ \forall t\in(t^n,t^{n+1}),\\
			\displaystyle\int_{\Omega} \frac{\partial u}{\partial t}(t) vdx +\int_{\Omega} \left[\left(\frac{\varepsilon}{2}\mathbf{I}+\cof(\p(t))\right)\nabla u(t)\right]\cdot \nabla v dx \, \\
			\quad\; +\, 2\int_{\Omega}K(1+|\s(t)|^2)^2 dx=0 \mbox{ in } \Omega\times(t^n,t^{n+1}),\; \forall v\in H^1_0(\Omega),\\
			u=g^n \mbox{ on } \partial \Omega\times (t^n,t^{n+1}) ,
		\end{cases}\\
		\frac{\partial \p}{\partial t}=\mathbf{0} \mbox{ in } \Omega\times(t^n,t^{n+1})  ,\\
		\frac{\partial \s}{\partial t}=\mathbf{0} \mbox{ in } \Omega\times(t^n,t^{n+1})  ,\\
		(u,\p,\s)(t^n)=(u^n,\p^n,\s^n),
	\end{cases}
	\label{eq.split2.1}
\end{equation}
and set
\begin{equation}
	u^{n+1/3}=u(t^{n+1})  ,\ \p^{n+1/3}=\p(t^{n+1}) ,\  \s^{n+1/3}=\s(t^{n+1}),\ g^{n+1/3}=g^n .
	\label{eq.split2.2}
\end{equation}
\underline{\emph{The Second Fractional Step:}}\\
Solve
\begin{equation}
	\begin{cases}
		\frac{\partial u}{\partial t}=0 \mbox{ in } \Omega\times(t^n,t^{n+1}) &\mbox{ in } \Omega\times(t^n,t^{n+1}),\\
		\frac{\partial \p}{\partial t}+\gamma_1(\p-\mathbf{D}^2 u^{n+1/3})=0 &\mbox{ in } \Omega\times(t^n,t^{n+1}),\\
		\frac{\partial \s}{\partial t}+\gamma_2(\s-\nabla u^{n+1/3})=0 &\mbox{ in } \Omega\times(t^n,t^{n+1}),\\
		(u,\p,\s)(t^n)=(u^{n+1/3} ,\p^{n+1/3},\s^{n+1/3}) ,
	\end{cases}
	\label{eq.split2.3}
\end{equation}
and set
\begin{equation}
	u^{n+2/3}=u(t^{n+1}) ,\  \p^{n+2/3}=P_+\left[\p(t^{n+1})\right],\ \s^{n+2/3}=\s(t^{n+1}), \ g^{n+2/3}=g^{n+1/3}.
	\label{eq.split2.4}
\end{equation}
\underline{\emph{The Third Fractional Step:}}\\
\begin{equation}
	\begin{cases}
		\begin{cases}
			u\in H^1(\Omega),\\
			\int_{\Omega} \frac{\partial u}{\partial t}(t)vdx + \frac{\varepsilon}{2} \int_{\Omega} \nabla u(t) \cdot \nabla vdx +\kappa\int_{\partial\Omega} (u(t)-g)vd\Gamma=0,\\
			\forall v\in H^1(\Omega),
		\end{cases}\\
		\frac{\partial \p}{\partial t}=\mathbf{0} \mbox{ in } \Omega\times(t^n,t^{n+1})  ,\\
		\frac{\partial \s}{\partial t}=\mathbf{0} \mbox{ in } \Omega\times(t^n,t^{n+1})  ,\\
		(u,\p,\s)(t^n)=(u^{n+2/3}  ,\p^{n+2/3} ,\s^{n+2/3}),
	\end{cases}
	\label{eq.split2.5}
\end{equation}
and set
\begin{equation}
	u^{n+1}=u(t^{n+1}) , \p^{n+1}=P_+\left[\p(t^{n+1})\right], \s^{n+1}=\s(t^{n+1}),g^{n+1}=u^{n+1}|_{\partial\Omega}.
	\label{eq.split2.6}
\end{equation}

Assuming that one uses just one step of the backward Euler scheme to solve the parabolic problem in (\ref{eq.split2.1}) and (\ref{eq.split2.5}), the Lie scheme (\ref{eq.split2.0})-(\ref{eq.split2.6}) reduces to the following variant of scheme (\ref{eq.split1Dis.0})-(\ref{eq.split1Dis.2}):
\begin{equation}
	(u^0,\p^0,\s^0,g^0)=(u_0,\p_0,\s_0,g).
	\label{eq.split2Dis.0}
\end{equation}
For $n\geq 0$, $(u^n,\p^n,\s^n,g^n)\rightarrow u^{n+1/2} \rightarrow (u^{n+1},\p^{n+1},\s^{n+1},g^{n+1})$ as follows:

\noindent Solve
\begin{eqnarray}
	&&\begin{cases}
		u^{n+1/2}\in H^1(\Omega),\\
		\displaystyle\int_{\Omega} \frac{u^{n+1/2}-u^n}{\Delta t}+\int_{\Omega} \left[\left(\frac{\varepsilon}{2}\mathbf{I}+
		\cof(\p^n)\right)\nabla u^{n+1/2}\right]\cdot \nabla vdx\,\\
		\quad \ +\,2\int_{\Omega}K(1+|\s^n|^2)^2vdx=0,\quad \forall v\in H^1_0(\Omega),\\
		u^{n+1/2}=g^n \mbox{ on } \partial\Omega,
	\end{cases}
	\label{eq.split2Dis.1}
\end{eqnarray}
and compute
\begin{eqnarray} 
	&&\begin{cases}
		\p^{n+1}=P_+\left[e^{-\gamma_1\Delta t}\p^n+
		\left(1-e^{-\gamma_1\Delta t}\right) \mathbf{D}^2 u^{n+1/2}\right], \\
		\s^{n+1}=e^{-\gamma_2\Delta t}\s^n+
		\left(1-e^{-\gamma_2\Delta t}\right) \nabla u^{n+1/2} ,\\
		\begin{cases}
			u\in H^1(\Omega),\\
			\displaystyle\int_{\Omega} \frac{u^{n+1}-u^{n+1/2}}{\Delta t}vdx + \varepsilon \int_{\Omega} \nabla u^{n+1} \cdot \nabla vdx\, \\
			\quad \ +\,\kappa\int_{\partial\Omega} (u^{n+1}-g^n)vd\Gamma=0, \quad \forall v\in H^1(\Omega),
		\end{cases}\\
		g^{n+1}=u^{n+1}|_{\partial\Omega}.
	\end{cases}
	\label{eq.split2Dis.2}
\end{eqnarray}

\section{Finite elements for the new operator-splitting scheme}
\label{sec.space}
The divergence form strongly suggests that we apply a finite-element method to implement (\ref{eq.split1Dis.1})-(\ref{eq.split1Dis.2}) and (\ref{eq.split2Dis.1})-(\ref{eq.split2Dis.2}). Here we choose a mixed finite-element method: we use the same function space to approximate $u$, $\nabla u$,
$\D^2u$, $\s$, and $\p$. Since we will choose basis functions to be piecewise affine functions, the resulting approximations are continuous piecewise affine on $\Omega$.

\subsection{Finite-element spaces}
Let $\mathcal{T}_h$ be the triangulation of the domain $\Omega$, and let $h$ denote the maximum edge length of the triangles
in $\mathcal{T}_h$. Let $\Sigma_h=\{Q_j\}_{j=1}^{N_h}$ be the collection of vertices in $\mathcal{T}_h$, where $Q_i$ denotes a typical vertex. We define the first finite-element space as
\begin{equation}
	V_h=\left\{v|v\in C^0\left(\bar{\Omega}\right),v|_T\in P_1,\forall T\in \mathcal{T}_h\right\},
	\label{eq.V}
\end{equation}
where $P_1$ denotes the space of polynomials with order no larger than $1$.

Accordingly, we associate each vertex $Q_j$ with a shape function $w_j$ such that
$$
w_j\in V_h, w_j(Q_j)=1, w_j(Q_k)=0,\;\;\forall k=1,\;\ldots,N_h,\; k\neq j,
$$
where the support of $w_j$, denoted $\omega_j$, is the union of triangles that have the same common vertex $Q_j$, and we denote the area of $\omega_j$ by $|\omega_j|$. The set $\mathcal{B}=\{w_j\}_{j=1}^{N_h}$ forms a collection of basis functions of $V_h$. In other words, we have
$$
v=\sum_{j=1}^{N_h} v(Q_j)w_j, \;\forall v\in V_h.
$$

In addition, we define
\begin{equation}
	V_{gh}=\left\{v|v\in V_h, v(Q_j)=g(Q_j),\;\forall Q_j\in \Sigma_h\cap\partial\Omega\right\},
	\label{eq.Vg}
\end{equation}
where $g$ can be any function which is $C^0$ on $\partial \Omega$. When $g=0$, we have that
$$
V_{0h}=V_h\cap H_0^1.
$$

Meanwhile, we define the following vector-valued spaces
\begin{eqnarray*}
	&&\mathbf{R}_h=\left\{\mathbf{r}|\;\mathbf{r}\in V_h^{2\times 1}\right\},\\
	&&\mathbf{Q}_h=\left\{\mathbf{q}|\;\mathbf{q}\in V_h^{2\times 2},\;\mathbf{q}=\mathbf{q}^T\right\},
\end{eqnarray*}
so that we can use functions in $\mathbf{R}_h$ to approximate $\nabla u$ and $\s$ and use functions in $\mathbf{Q}_h$ to approximate $\D^2u$ and $\p$.

\subsection{Approximations of the two first-order derivatives of $u$}
\label{sec.gradApprox}
For any $v\in V_h$, we denote the first-order derivative approximation $\frac{\partial v}{\partial x_i}$ of $v$ by $D_{ih}(v)$ for $i=1,2$, and this approximate derivative operator is defined in the following weak sense:
\begin{eqnarray}
	\int_{\Omega} D_{ih}(v) wdx=\int_{\Omega} \frac{\partial v}{\partial x_i} w dx, \;\; i=1,2, \;\forall w\in H^1(\Omega).
	\label{eq.gradApp}
\end{eqnarray}
Since $\Omega$ is partitioned by the triangulation $\mathcal{T}_h$, we restrict the test functions $w$ to be in $V_h$ so that we only need to test the above integral against those basis functions $w_k$ for $k=1,2, \cdots, N_h$. Since $w_k$ is only supported on $\omega_k$, we have
\begin{equation}
	\begin{cases}
		D_{ih}(v)\in V_h,\;\;\forall i=1,2,\\
		D_{ih}(v)(Q_k)=\frac{3}{|\omega_k|}\int_{\omega_k} \frac{\partial v}{\partial x_i}w_kdx,\;\;\forall k=1,2,...,N_h.
	\end{cases}
	\label{eq.gradApprox}
\end{equation}

We remark in passing that on a regular mesh such as the one shown in Figure \ref{fig.4mesh}(a), (\ref{eq.gradApprox}) recovers the central-difference approximation at an interior node and one-sided approximation at a boundary node in a finite-difference method based on this mesh.

In some problems, $\nabla u$ has singularities on $\Omega$. One challenging situation is when the singularities appear on the boundary. Since the above approximation is second-order accurate at interior nodes and first-order accurate at boundary nodes, the approximation at nodes near the boundary can be very unstable and even blow up, especially when the gradient of the exact solution blows up at the boundary of a computational domain, such as a semi-sphere. To resolve this problem, we need to regularize the approximation of $\nabla u$. One possible way is to adopt the idea from \cite{glowinskillq2017,caboussat2013least} which is used to approximate the second-order derivative:
\begin{equation}
	\begin{cases}
		D_{ih}(v)\in H_0^1,\\
		\epsilon \int_{\Omega} \nabla D_{ih}\cdot \nabla wdx+\int_{\Omega} D_{ih}(v) wdx=\int_{\Omega} \frac{\partial v}{\partial x_i} w dx, \;\; i=1,2, \;\forall w\in H^1_0(\Omega).
	\end{cases}
	\label{eq.firstapprox0}
\end{equation}
The error of the regularized approximation can be larger than that of the direct approximation, but it is more robust. Moreover, we have
$$
\lim_{\epsilon,h\rightarrow0} D_{ih}(v)=\frac{\partial v}{\partial x_i} \;\;\mbox{ in } L^2(\Omega).
$$

\subsection{Approximations of second-order derivatives of $u$}
The general idea to approximate the second-order derivatives is similar to the one used in \cite{glowinskillq2017}. For completeness, we mention the details here.

For any $v\in V_h$, we denote the approximations of $\frac{\partial^2 v}{\partial x_i \partial x_j}$ by $D^2_{ijh}(v)$ for $i,j=1,2$, so that the approximate operator $D^2_{ijh}(v)$ of second-order derivatives $D^2_{ijh}(v)$ is defined in the following weak sense,
\begin{equation}
	\int_{\Omega} D^2_{ijh}(v)w_k dx=\int_{\Omega} \frac{\partial^2 v}{\partial x_i \partial x_j}w_k dx.
	\label{eq.secondApprox0}
\end{equation}

To resolve the right hand side of (\ref{eq.secondApprox0}), we apply the divergence theorem,
\begin{equation}
	\int_{\Omega} \frac{\partial^2 v}{\partial x_i \partial x_j}w dx=\frac{1}{2}\int_{\partial \Omega} \left(\frac{\partial v}{\partial x_i}n_j + \frac{\partial v}{\partial x_j} n_i\right)w d(\partial\Omega)-\frac{1}{2}\int_{\Omega} \left(\frac{\partial v}{\partial x_i} \frac{\partial w}{\partial x_j} + \frac{\partial v}{\partial x_j} \frac{\partial w}{\partial x_i}\right) dx,
	\label{eq.secondApprox1}
\end{equation}
where $\mathbf{n}=(n_1,n_2)$ is the outward normal direction along $\partial\Omega$. The above approximation is accurate at interior nodes, but the approximation error is large at nodes on the boundary. For example, consider the approximate derivative operator $D^2_{11h}$ on a regular mesh of the unit square; after some derivation, we can show that there is always one node at one of the corners of the unit square such that $D^2_{11h}(v)=0$ at that node, no matter what form $v$ takes.

To deal with this issue, we treat interior nodes and boundary nodes separately. Let $\Sigma_{0h}=\{Q_k\}_{k=1}^{N_0}$ denote the set of interior nodes in $\Omega$, where we assume that the first $N_0$ nodes of $\Sigma_h$ are in the interior of $\Omega$. It follows that we have $\Sigma_h\cap\partial\Omega=\{Q_k\}_{k=N_0+1}^{N_h}$. For $k=1,2,...,N_0$, the approximation of (\ref{eq.secondApprox0})-(\ref{eq.secondApprox1}) reduces to
\begin{eqnarray}
	\int_{\Omega} D^2_{ijh}(v)w_k dx=-\frac{1}{2}\int_{\Omega} \left(\frac{\partial v}{\partial x_i} \frac{\partial w_k}{\partial x_j} + \frac{\partial v}{\partial x_j} \frac{\partial w_k}{\partial x_i}\right) dx.
	\label{eq.crime1}
\end{eqnarray}

To treat nodes on the boundary, the work in \cite{caboussat2013least} used the zero Dirichlet boundary condition for the operator  $D^2_{ijh},i,j=1,2$, though the boundary value is not needed in the resulting algorithm. In comparison with the numerical method  in \cite{caboussat2013least}, ours are different in that the boundary value is crucial for our splitting algorithm. Specifically, in (\ref{eq.split1Dis.1})-(\ref{eq.split1Dis.2}) we need boundary values to update $\p$ which is in turn used to update $u$. Therefore, we need a better treatment of the boundary nodes.

Here we adopt a strategy from \cite{glowinskillq2017,glowinskillq20173d}  to treat boundary nodes by committing a ``variational crime''. First, we impose the zero Neumann boundary condition
\begin{equation}
	\frac{\partial D^2_{ijh}(v)}{\partial \mathbf{n}}=0.
	\label{eq.neu}
\end{equation}
Multiplying (\ref{eq.neu}) by $w_k$ for $k=N_0+1,...,N_h$ and integrating along $\partial\Omega$, we get
\begin{eqnarray}
	0&=&\int_{\partial\Omega}\frac{\partial D^2_{ijh}(v)}{\partial \mathbf{n}}w_k d(\partial\Omega)=\int_{\Omega} \nabla\cdot \left( \nabla D^2_{ijh}(v) w_k\right)dx\nonumber\\
	&=&\int_{\Omega} \nabla^2 D^2_{ijh}(v) w_kdx +\int_{\Omega} \nabla D^2_{ijh}(v)\cdot \nabla w_kdx.
	\label{eq.neu1}
\end{eqnarray}
If $D^2_{ijh}(v)$ is harmonic, implying that $\nabla^2 D^2_{ijh}(v)=0$, then we have
\begin{equation}
	\int_{\Omega} \nabla D^2_{ijh}(v)\cdot \nabla w_kdx=0.
	\label{eq.crime2}
\end{equation}
In our algorithm, although $D^2_{ijh}$ is only piecewise harmonic, we still use (\ref{eq.crime2}) to update boundary values, which is the so-called variational crime. In either approximation (\ref{eq.secondApprox0})-(\ref{eq.secondApprox1}) or approximation   (\ref{eq.crime1}) and (\ref{eq.crime2}), since $w_k$ is only supported on $\omega_k$, the integration domain can be replaced by $\omega_k$ if the test function is $w_k$.

In our numerical experiments, with the regularization mechanism introduced below, the accuracy by (\ref{eq.crime1}) and (\ref{eq.crime2}) is similar to that by (\ref{eq.secondApprox0})-(\ref{eq.secondApprox1}), but (\ref{eq.crime1}) and (\ref{eq.crime2}) make the algorithm more robust. It is worth mentioning that as implemented in \cite{glowinskillq2017} both approximations work  for two-dimensional Monge-Amp\`{e}re equations; however, as shown in \cite{glowinskillq20173d} only the approximation based on the variational crime works for three-dimensional Monge-Amp\`{e}re equations.

As reported in \cite{caboussat2013least,glowinskillq2017,glowinskillq20173d}, if we directly use the above approximations, the performance of our algorithm depends on triangulations; in the worst case, on a symmetric mesh as shown in Figure \ref{fig.4mesh}(b), our algorithm does not converge. To obtain an algorithm which is robust for all kinds of meshes, we need to regularize
the problem by adding some viscosity to our formulation of second-order derivatives.

As a first approach of regularization, we incorporate a local regularization term into the weak definition of second-order derivatives at interior nodes:
\begin{equation}
	\begin{cases}
		\forall i,j=1,2,\;\forall v\in V_h, \; D^2_{ijh}(v)\in V_h \mbox{ and }\\
		C\sum_{T\in\mathcal{T}^k_h}|T|\int_T \nabla D^2_{ijh}(v)\cdot\nabla w_k dx+\int_{\omega_k}D^2_{ijh}(v) w_k dx\\
		\hspace{1cm}=-\frac{1}{2}\int_{\omega_k}\left[ \frac{\partial v}{\partial x_i} \frac{\partial w_k}{\partial x_j}+\frac{\partial v}{\partial x_j} \frac{\partial w_k}{\partial x_i}\right] dx, \; \forall k=1,...,N_{0h},\\
		\int_{\omega_k} \nabla D^2_{ijh}(v)\cdot\nabla w_k dx=0, \; \forall k=N_{0h}+1,...,N_h,
	\end{cases}
	\label{eq.smooth}
\end{equation}
and the above relations at the discretization level can be written as
\begin{equation}
	\begin{cases}
		\forall i,j=1,2,\; \forall v\in V_h, \; D^2_{ijh}(v)\in V_h \mbox{ and }\\
		C\sum_{T\in\mathcal{T}^k_h}|T|\int_T \nabla D^2_{ijh}(v)\cdot\nabla w_k dx+\int_{\omega_k}D^2_{ijh}(v) w_k dx\\
		\hspace{1cm}=-\frac{1}{2}\int_{\omega_k}\left[ \frac{\partial v}{\partial x_i} \frac{\partial w_k}{\partial x_j}+\frac{\partial v}{\partial x_j} \frac{\partial w_k}{\partial x_i}\right] dx, \; \forall k=1,...,N_{0h},\\
		\int_{\omega_k} \nabla D^2_{ijh}(v)\cdot\nabla w_k dx=0, \;\forall k=N_{0h}+1,...,N_h,
	\end{cases}
	\label{eq.smoothDis}
\end{equation}
where $C$ is a positive constant of order $1$, $\mathcal{T}_h^k$ is the set of all triangles with the common vertex $Q_k$.

If all triangles in $\mathcal{T}_h$ are of a similar size, (\ref{eq.smooth}) can be slightly simplified to be
\begin{equation}
	\begin{cases}
		\forall i,j=1,2,\; \forall v\in V_h, D^2_{ijh}(v)\in V_h \mbox{ and }\\
		\epsilon_2\int_{\omega_k} \nabla D^2_{ijh}(v)\cdot\nabla w_k dx+\int_{\omega_k}D^2_{ijh}(v) w_k dx\\
		\hspace{1cm}=-\frac{1}{2}\int_{\omega_k}\left[ \frac{\partial v}{\partial x_i} \frac{\partial w_k}{\partial x_j}+\frac{\partial v}{\partial x_j} \frac{\partial w_k}{\partial x_i}\right] dx, \;\forall k=1,...,N_{0h},\\
		\int_{\omega_k} \nabla D^2_{ijh}(v)\cdot\nabla w_k dx=0, \;\forall k=N_{0h}+1,...,N_h,
	\end{cases}
	\label{eq.smoothSimp}
\end{equation}
where $\varepsilon_2$ is of order $O(h^2)$.

As a second approach of regularization, we incorporate a double-regularization mechanism into our weak formulation of second-order derivatives. Assuming that $\psi\in H^2$, we consider the following linear elliptic variational problem,
\begin{equation}
	\begin{cases}
		p_{ij}^{\varepsilon}\in H_0^1(\Omega),\\
		\varepsilon\int_{\Omega} \nabla p_{ij}^{\varepsilon} \cdot \nabla \phi dx+ \int_{\Omega} p_{ij}^{\varepsilon}\phi dx= -\frac{1}{2} \displaystyle\int_{\Omega} \left[ \frac{\partial\psi}{\partial x_i}\frac{\partial \phi}{\partial x_j} + \frac{\partial \psi}{\partial x_j} \frac{\partial \phi}{\partial x_i} \right]  dx, \forall \phi\in H_0^1(\Omega),
	\end{cases}
	\label{eq.double.1}
\end{equation}
which yields the following relations in the weak sense
\begin{equation}
	\lim_{\varepsilon\rightarrow0} p_{ij}=\frac{\partial^2 \psi}{\partial x_i \partial x_j} \;\;\mbox{ in } L^2(\Omega),
	\label{eq.double.2}
\end{equation}
and
\begin{equation}
	\begin{cases}
		-\varepsilon\nabla^2 p_{ij}^{\varepsilon}+ p_{ij}^{\varepsilon} = \frac{\partial^2 \psi}{\partial x_i \partial x_j}\;\;\mbox{ in } \Omega,\\
		p_{ij}^{\varepsilon}=0 \;\;\mbox{ on } \partial \Omega.
	\end{cases}
	\label{eq.double.3}
\end{equation}
Since, as reported in \cite{glowinskillq2017}, this approximation is not effective in treating the zero-Dirichlet boundary condition,  we apply the following correction step,
\begin{equation}
	\begin{cases}
		-\varepsilon\nabla^2 \tilde{p}_{ij}^{\varepsilon}+ \tilde{p}_{ij}^{\varepsilon}= p_{ij}^{\varepsilon} \;\;\mbox{ in } \Omega,\\
		\frac{\partial \tilde{p}_{ij}^{\varepsilon}}{\partial \n}=0 \;\; \mbox{ on } \partial\Omega,
	\end{cases}
	\label{eq.double.4}
\end{equation}
whose variational formulation reads as
\begin{equation}
	\begin{cases}
		\tilde{p}_{ij}^{\varepsilon}\in H^1(\Omega),\\
		\varepsilon\int_{\Omega} \nabla \tilde{p}_{ij}^{\varepsilon} \cdot \nabla \phi dx + \int_{\Omega} \tilde{p}_{ij}^{\varepsilon} \phi dx=\int_{\Omega} p_{ij}^{\varepsilon} \phi dx, \;\; \forall \phi \in H^1(\Omega).
	\end{cases}
	\label{eq.double.5}
\end{equation}
It follows that $\tilde{p}_{ij}^{\varepsilon}$ verifies $\lim_{\varepsilon\rightarrow0} \tilde{p}_{ij}^{\varepsilon}= \frac{\partial^2 \psi}{\partial x_i \partial x_j}$ in $L^2(\Omega)$, and $\tilde{p}_{ij}^{\varepsilon} \in H^4(\Omega)$.

Consequently, the discrete analogue $D^2_{ijh}(v)$ of $\frac{\partial^2 v}{\partial x_i \partial x_j} (1\leq i,j \leq 2)$ can be computed in the following way: \\
Solve:
\begin{equation}
	\begin{cases}
		p_{ij}\in V_{0h},\\
		C\sum_{T\in \mathcal{T}_h^k} |T| \int_T \nabla p_{ij} \cdot \nabla w_k dx +\frac{|\omega_k|}{3}p_{ij}(Q_k)= -\frac{1}{2}\int_{\omega_k} \left[ \frac{\partial v}{\partial x_i} \frac{\partial w_k}{\partial x_j} + \frac{\partial v}{\partial x_j} \frac{\partial w_k}{\partial x_i}\right]dx,\\
		\forall k=1,...,N_{0h},
	\end{cases}
	\label{eq.double.6}
\end{equation}
and then
\begin{equation}
	\begin{cases}
		D_{ijh}^2(v) \in V_h,\\
		C\sum_{T\in \mathcal{T}_h^k} |T| \int_T \nabla D_{ijh}^2(v) \cdot \nabla w_k dx +\frac{|\omega_k|}{3}D_{ijh}^2(v)(Q_k)= \frac{|\omega_k|}{3}p_{ij}(Q_k),\\
		\forall k=1,...,N_{h},
	\end{cases}
	\label{eq.double.7}
\end{equation}
where $C$ is a constant of order $1$. Similar to the first approach of regularization, if all triangles in $\mathcal{T}_h$ are of a similar size, we can replace $C\sum_{T\in \mathcal{T}_h^k} |T|$ in (\ref{eq.double.6}) and (\ref{eq.double.7}) by $\epsilon_2$ which is of order $O(h^2)$.

\subsection{Implementation of scheme (\ref{eq.split1Dis.0})-(\ref{eq.split1Dis.2})}
We give a fully discretized analogue of scheme (\ref{eq.split1Dis.0})-(\ref{eq.split1Dis.2}) as follows. \\
Initialize
\begin{equation}
	u^0=u_0\in V_h,\p^0=\p_0\in\mathbf{Q}_h, \s^0=\s_0\in \mathbf{R}_h.
	\label{eq.split1Sch.0}
\end{equation}
For $n\geq 0$, proceed $\{u^n,\p^n,\s^n\}\rightarrow\{u^{n+1},\p^{n+1},\s^{n+1}\}$ as the following.  \\
Solve
\begin{equation}
	\begin{cases}
		u^{n+1}\in V_{gh},\\
		\int_{\Omega} u^{n+1}vdx+
		\Delta t\int_{\Omega}\left(\varepsilon \mathbf{I}+\cof(\p^n)\right)\nabla u^{n+1} \cdot \nabla vdx \\
		\hspace{1cm} = \int_{\Omega} u^n vdx-2\Delta tK\int_{\Omega}(1+|\s^n|^2)^2dx, \forall v\in V_{0h},
	\end{cases}
	\label{eq.split1Sch.1}
\end{equation}
and compute $\p^{n+1}$ and $\s^{n+1}$ via
\begin{equation}
	\begin{cases}
		\forall k=1,...,N_h,\\
		\alpha = e^{-\gamma_1\Delta t},	\\	
		\p^{n+\frac{1}{2}}(Q_k)=\alpha\; \p^n(Q_k)\,+\,( 1-\alpha)
		\begin{pmatrix}
			D_{11h}^2( u^{n+1})(Q_k) & D_{12h}^2( u^{n+1})(Q_k)\\
			D_{12h}^2( u^{n+1})(Q_k) & D_{22h}^2( u^{n+1})(Q_k)
		\end{pmatrix},\\
		\p^{n+1}(Q_k)=P_+\left[ \p^{n+1/2}(Q_k)\right].
	\end{cases}
	\label{eq.split1Sch.2}
\end{equation}
and
\begin{equation}
	\begin{cases}
		\forall k=1,...,N_h,\\
		\s^{n+1}(Q_k)=e^{-\gamma_2\Delta t}\s^n(Q_k)+\left( 1-e^{-\gamma_2 \Delta t}\right)
		\begin{pmatrix}
			D_{1h}(u^{n+1}(Q_k)\\
			D_{2h}(u^{n+1}(Q_k)
		\end{pmatrix}.
	\end{cases}
	\label{eq.split1Sch.3}
\end{equation}
Here, all integrations in (\ref{eq.split1Sch.1}) are computed by the trapezoidal rule. In (\ref{eq.split1Sch.2}) and (\ref{eq.split1Sch.3}), $D^1_{ih}$ for $i=1,2$ are computed using (\ref{eq.firstapprox0}) or (\ref{eq.gradApp}); $D^2_{ijh}(u^{n+1})$ for $i,j=1,2$ are  computed by approximation (\ref{eq.smooth}) or (\ref{eq.double.6})-(\ref{eq.double.7}).

\subsection{Implementation of scheme (\ref{eq.split2Dis.0})-(\ref{eq.split2Dis.2})}
The discretized analogue of scheme (\ref{eq.split1Dis.0})-(\ref{eq.split1Dis.2}) can be written as: \\
Initialize
\begin{equation}
	u^0=u_0\in V_h,\, \p^0=\p_0\in\mathbf{Q}_h,\, \s^0=\s_0\in \mathbf{R}_h, \, g^0=g.
	\label{eq.split2Sch.0}
\end{equation}
For $n\geq 0$, proceed $\{u^n,\p^n,\s^n\}\rightarrow\{u^{n+1},\p^{n+1},\s^{n+1}\}$ as the following.  \\
Solve
\begin{equation}
	\begin{cases}
		u^{n+1/2}\in V_{g^nh},\\
		\int_{\Omega} u^{n+1/2}vdx+
		\Delta t\int_{\Omega}\left(\varepsilon \mathbf{I}+\cof(\p^n)\right)\nabla u^{n+1/2} \cdot \nabla vdx \\
		\hspace{1cm} = \int_{\Omega} u^n vdx-2\Delta tK\int_{\Omega}(1+|\s^n|^2)^2dx, \forall v\in V_{g^nh}.
	\end{cases}
	\label{eq.split2Sch.1}
\end{equation}
Compute $\p^{n+1}$ and $\s^{n+1}$ via
\begin{equation}
	\begin{cases}
		\forall k=1,...,N_h,\\
		\alpha = e^{-\gamma_1\Delta t}, \\		
		\p^{n+\frac{1}{2}}(Q_k)=\alpha\, \p^n(Q_k)\, +\, ( 1-\alpha)
		\begin{pmatrix}
			D_{11h}^2\left( u^{n+1/2}\right)(Q_k) & D_{12h}^2\left( u^{n+1/2}\right)(Q_k)\\
			D_{12h}^2\left( u^{n+1/2}\right)(Q_k) & D_{22h}^2\left( u^{n+1/2}\right)(Q_k)
		\end{pmatrix},\\
		\p^{n+1}(Q_k)=P_+\left[ \p^{n+1/2}(Q_k)\right].
	\end{cases}
	\label{eq.split2Sch.2}
\end{equation}
and
\begin{equation}
	\begin{cases}
		\forall k=1,...,N_h,\\
		\s^{n+1}(Q_k)=e^{-\gamma_2\Delta t}\s^n(Q_k)+\left( 1-e^{-\gamma_2 \Delta t}\right)
		\begin{pmatrix}
			D_{1h}(u^{n+1/2}(Q_k)\\
			D_{2h}(u^{n+1/2}(Q_k)
		\end{pmatrix}.
	\end{cases}
	\label{eq.split2Sch.3}
\end{equation}
Compute
\begin{equation}
	\begin{cases}
		u^{n+1}\in V_{h},\\
		\int_{\Omega} u^{n+1}vdx+
		\Delta t\varepsilon\int_{\Omega} \nabla u^{n+1} \cdot \nabla vdx +\Delta t \int_{\partial\Omega} u^{n+1}vdx\\
		\hspace{1cm} = \int_{\Omega} u^{n+1/2} vdx + \Delta t \varepsilon\int_{\partial\Omega} gvdx, \forall v\in V_{h}.
	\end{cases}
	\label{eq.split2Sch.4}
\end{equation}
and update
\begin{equation}
	g^{n+1}=u^{n+1}|_{\partial\Omega}.
	\label{eq.split2Sch.5}
\end{equation}
All integrations in (\ref{eq.split2Sch.1}) and (\ref{eq.split2Sch.4}) are computed by the trapezoidal rule. In (\ref{eq.split2Sch.2}) and (\ref{eq.split2Sch.3}), $D^1_{ih}$ for $i=1,2$ are computed using (\ref{eq.firstapprox0}) or (\ref{eq.gradApp}); $D^2_{ijh}(u^{n+1})$ for $i,j=1,2$ are  computed by approximation (\ref{eq.smooth}) or (\ref{eq.double.6})-(\ref{eq.double.7}).

\subsection{The projection operator $P_+(\cdot)$}
\label{sec.P}
Since we want to find a convex solution $u$, we need to have some mechanism to enforce convexity in our algorithm. There are many possible approaches to handle the issue.

One particular approach that we discuss here is to modify one of the finite-element components, $\p$, after each iteration so that the modified $\p$ satisfies some convexity-related properties. Since the Hessian matrix of a convex function is semi-positive definite and we expect $\p$ to converge to the Hessian matrix of the exact solution $u^*$ which is convex, it is reasonable to enforce $\p$ to be semi-positive definite; therefore, we introduce a spectral projection operator to achieve this, and $P_+(\cdot)$ is such a projector in our algorithm.

Let $\mathbf{A}$ be a symmetric $2\times2$ matrix. Assume that $\mathbf{A}$ has a spectral decomposition, $\mathbf{A}=\mathbf{S}\mathbf{\Lambda}\mathbf{S}^{-1}$, where the columns of $\mathbf{S}$ are the eigenvectors of $\mathbf{A}$ and $\mathbf{\Lambda}=
\begin{pmatrix}
	\lambda_1 &0\\
	0 & \lambda_2
\end{pmatrix}.$
We define the spectral projector operator  $P_+(\cdot)$ as
$$
P_+(\mathbf{A})=\mathbf{S}
\begin{pmatrix}
	\lambda_1^+ &0\\
	0 & \lambda_2^+
\end{pmatrix}
\mathbf{S}^{-1},
$$
where $\lambda_i^+=\max\{\lambda_i,0\}$ for $i=1,2$. The effect of $P_+(A)$ is to project $A$ onto the cone consisting of semi-positive definite matrices. This projection during each iteration makes equation (\ref{eq.split1Dis.1}) an elliptic PDE of $u$.

Another possible approach is to choose a convex initial condition which will be discussed in the next section.

\section{Initialization}
\label{sec.initial}
\subsection{Initial condition for scheme (\ref{eq.split1Dis.0})-(\ref{eq.split1Dis.2})}
To initialize $u_0$ and $\p_0$ for scheme (\ref{eq.split1Dis.0})-(\ref{eq.split1Dis.2}), we solve the standard Monge-Amp\`{e}re equation
\begin{equation}
	\begin{cases}
		\det(\D^2u_0)=K, \\
		u_0=g \;\mbox{ on } \partial\Omega.
	\end{cases}
	\label{eq.initial.u}
\end{equation}
We will deal with (\ref{eq.initial.u}) by adopting the method in \cite{glowinskillq2017}, which solves the following initial value problem to steady state,
\begin{equation}
	\begin{cases}
		\begin{cases}
			\frac{\partial u}{\partial t}-\nabla\cdot\left((\varepsilon\mathbf{I}+\cof(\p))\nabla u\right)+2K=0,\\
			u=g \mbox{ on } \partial\Omega,
		\end{cases}\\
		\frac{\partial \p}{\partial t}+\gamma(\p-\D^2u)=\mathbf{0}.
	\end{cases}
	\label{eq.ma2d}
\end{equation}

Let $\{u_*,\p_*\}$ be the steady state of (\ref{eq.ma2d}). Accordingly, we set $u_0=u_*$, $\p_0=\D^2u_*$ and $\s=\D u_*$ as the initial condition for our scheme (\ref{eq.split1Dis.0})-(\ref{eq.split1Dis.2}). Therefore, our algorithm can be summarized as a two-stage method:
\begin{center}
	\parbox{0.9\textwidth}{
		\textbf{\emph{Stage 1}}\\
		In the algorithm in \cite{glowinskillq2017}, set $\varepsilon_1=\varepsilon_2=h^2$ and $dt=2h^2$. Solve (\ref{eq.ma2d}) until $\|u^{n+1}-u^n\|_2<tol$ to get $u_0$. Compute $\p_0=\D^2u_0$ and $\s_0=\D u_0$.\\
		\textbf{\emph{Stage 2}}\\
		With the initial condition $u_0$, $\p_0$, and $\s_0$, solve (\ref{eq.split1Dis.0})-(\ref{eq.split1Dis.2}) to steady state.
	}
\end{center}

\subsection{Initial condition for scheme (\ref{eq.split2Dis.0})-(\ref{eq.split1Dis.2})}
When using scheme (\ref{eq.split2Dis.0})-(\ref{eq.split1Dis.2}), the boundary value of the computed solution does not satisfy the given boundary condition, so the initial condition used for scheme (\ref{eq.split1Dis.0})-(\ref{eq.split1Dis.2}) may not help. To initialize scheme (\ref{eq.split2Dis.0})-(\ref{eq.split1Dis.2}), the initial condition used to solve (\ref{eq.ma2d}) in \cite{glowinskillq2017}:
\begin{equation}
	\begin{cases}
		\nabla^2 u_0=2\lambda\sqrt{K},\\
		u_0|_{\partial\Omega}=g,
	\end{cases}
	\label{eq.ma2d.initial}
\end{equation}
where $\lambda$ $(>0)$ is of order $O(1)$.

\begin{figure}[!ht]
	\centering
	\begin{tabular}{cccc}
		(a) & (b) & (c) &(d)\\
		\includegraphics[width=0.15\textwidth]{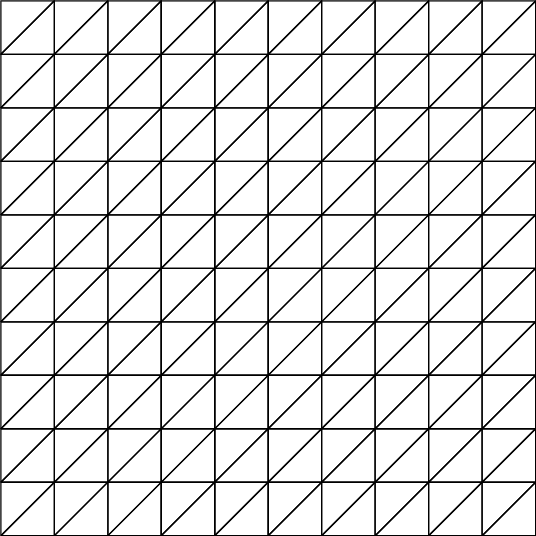} &
		\includegraphics[width=0.15\textwidth]{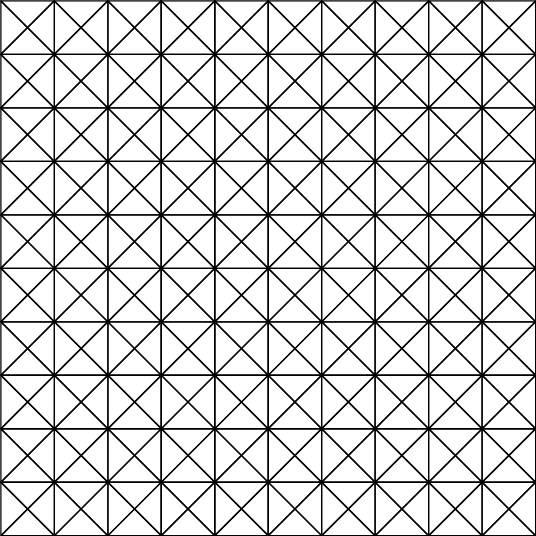} &
		\includegraphics[width=0.15\textwidth]{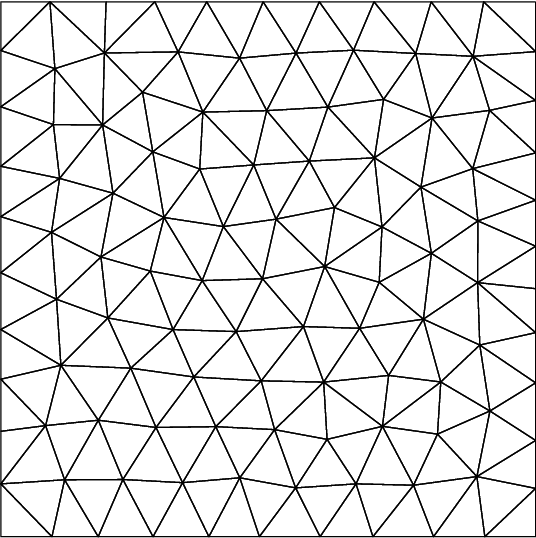}&
		\includegraphics[width=0.15\textwidth]{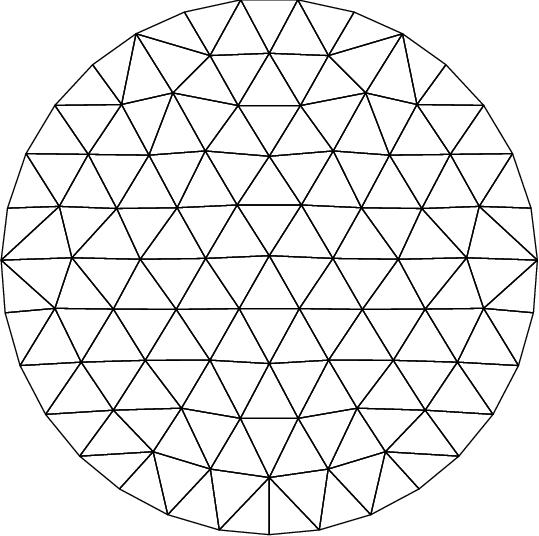}
	\end{tabular}
	
	\caption{Four meshes for two different domains used in our numerical experiments. (a) A regular mesh on a square. (b) A (highly) symmetric mesh on a square. (c) An anisotropic unstructured mesh on a square. (d) An anisotropic unstructured mesh on a half-unit disk.}
	\label{fig.4mesh}
\end{figure}

\section{Numerical experiments}
\label{sec.example}
In this section, we carry out a variety of numerical experiments in different settings to demonstrate the performance of scheme (\ref{eq.split1Dis.0})-(\ref{eq.split1Dis.2}) and scheme (\ref{eq.split2Dis.0})-(\ref{eq.split2Dis.2}). Four different meshes as shown in Figure \ref{fig.4mesh} will be used in our experiments: (a) regular meshes on a unit square, (b) symmetric meshes on a unit square, (c) unstructured meshes on a unit square, and (d) unstructured meshes on a half-unit disk. In all of our experiments, in Stage 1 of our algorithm, we use the method in \cite{glowinskillq2017} to initialize the iteration of our algorithm, where we use $tol=h^2$. Without specification, we choose $\Delta t=2h^2$ and $\varepsilon_1=\varepsilon_2=h^2$ in both Stage 1 and 2 of our algorithm. For examples with compatible boundary condition, scheme (\ref{eq.split1Dis.0})-(\ref{eq.split1Dis.2}) is used. For examples with incompatible boundary condition, scheme (\ref{eq.split2Dis.0})-(\ref{eq.split2Dis.2}) is used. We also compare the numerical solutions by both schemes on some examples. Without specification, stopping criterion $\|u^{n+1}-u^n\|_2<10^{-6}$ and scheme (\ref{eq.split1Dis.0})-(\ref{eq.split1Dis.2}) are used. 

\begin{figure}
	\centering
	\begin{tabular}{cccc}
		(a) & (b) & (c) & (d)\\
		\includegraphics[width=0.21\textwidth]{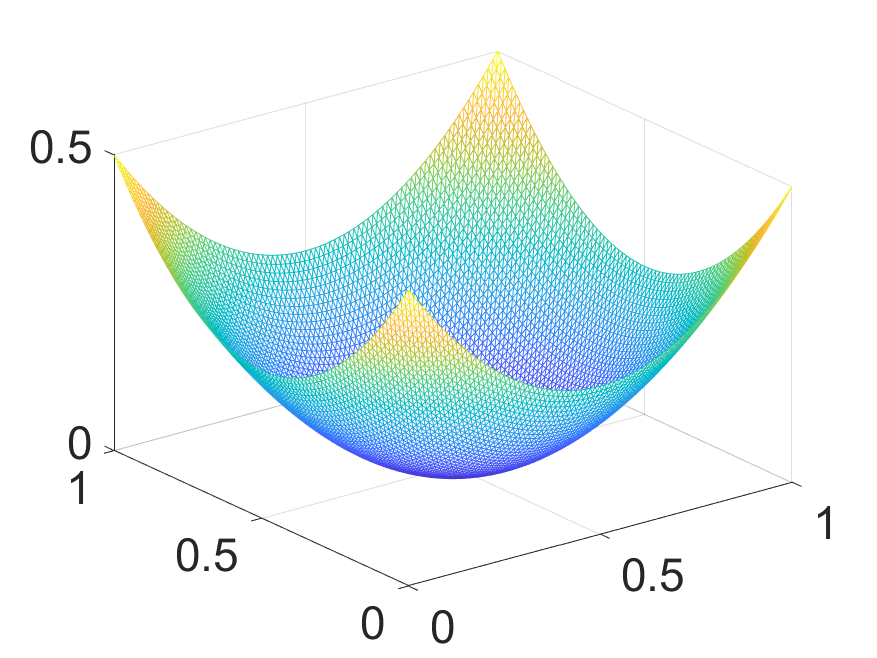} & 
		\includegraphics[width=0.21\textwidth]{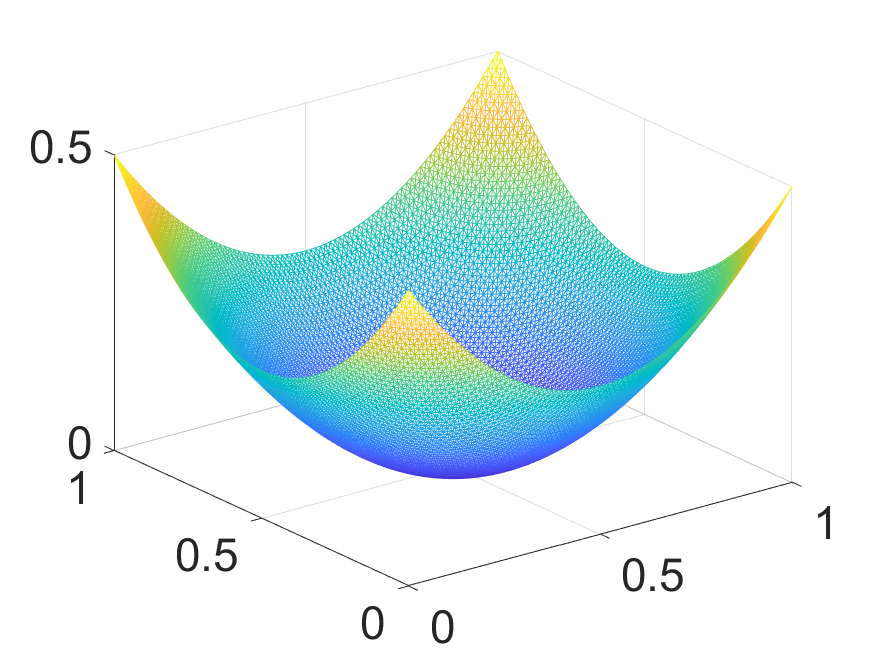}&
		\includegraphics[width=0.21\textwidth]{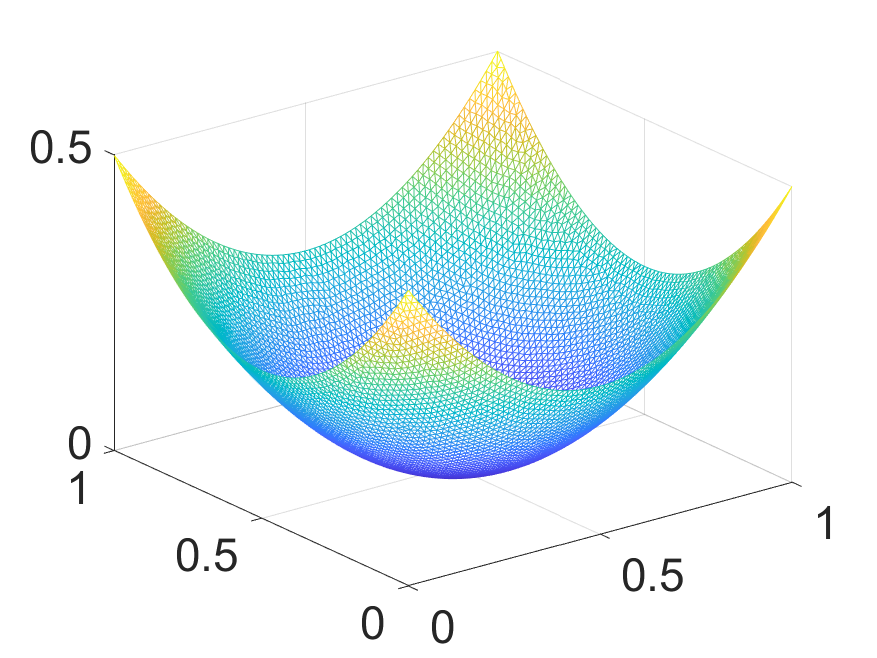}&
		\includegraphics[width=0.21\textwidth]{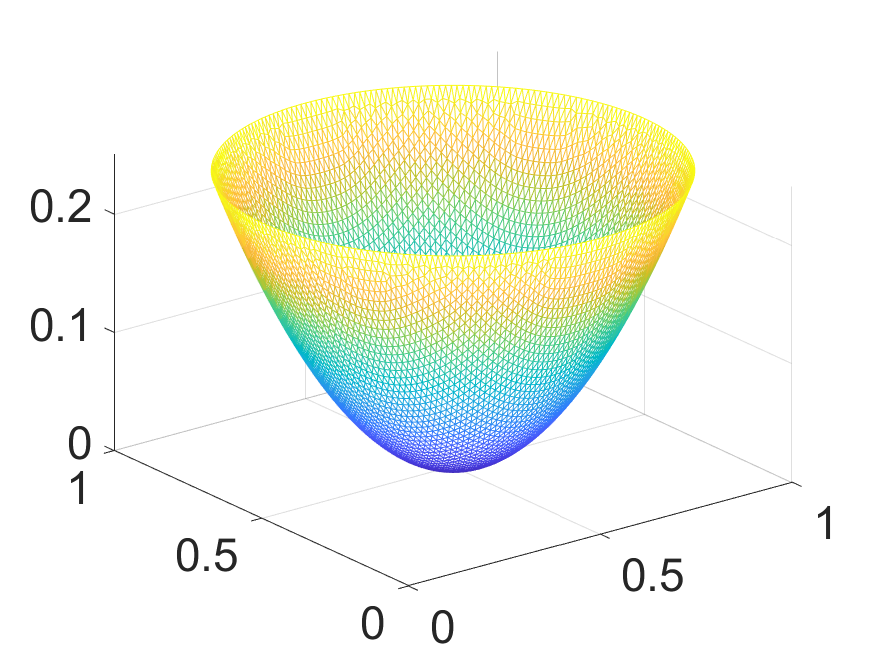}\\
		\includegraphics[width=0.21\textwidth]{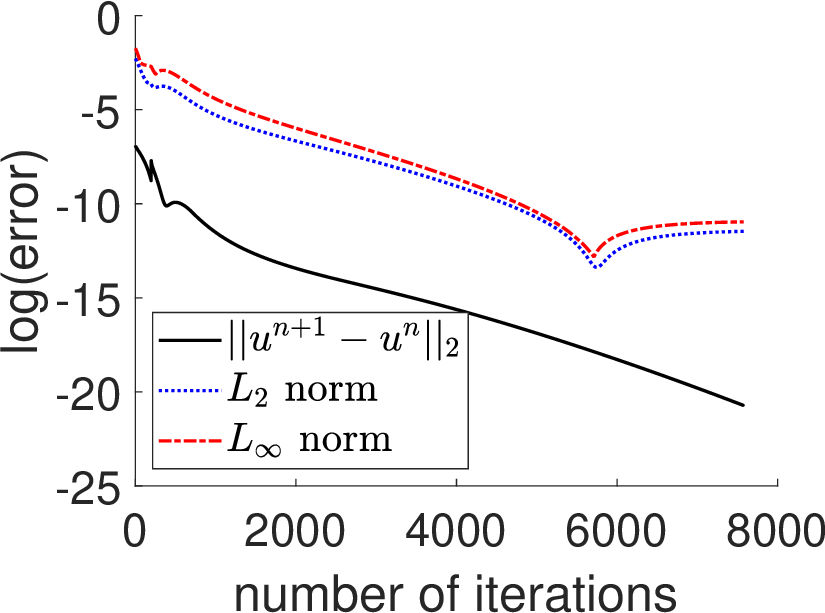}&
		\includegraphics[width=0.21\textwidth]{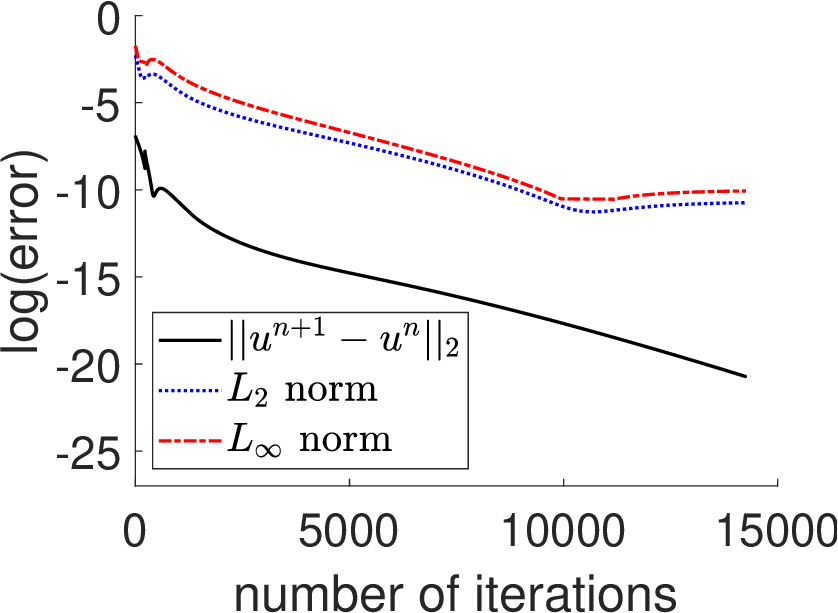}&
		\includegraphics[width=0.21\textwidth]{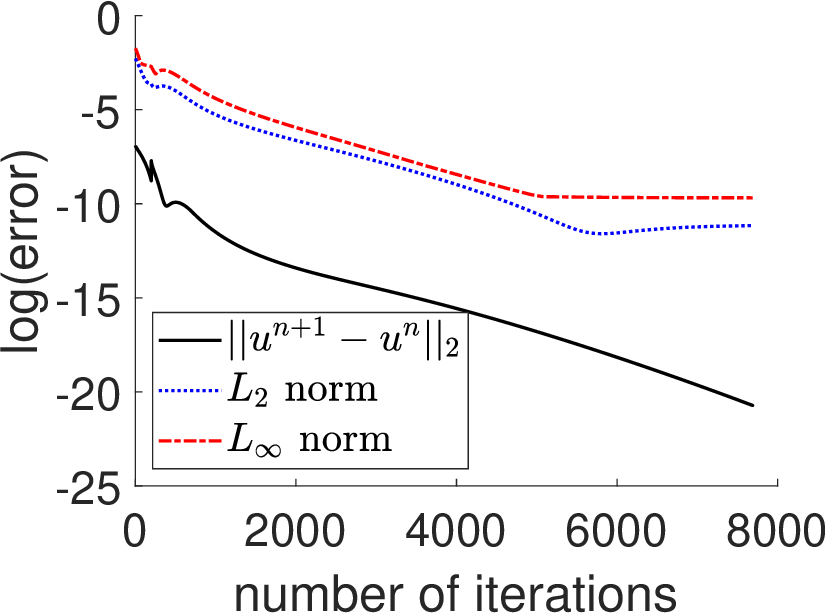}&
		\includegraphics[width=0.21\textwidth]{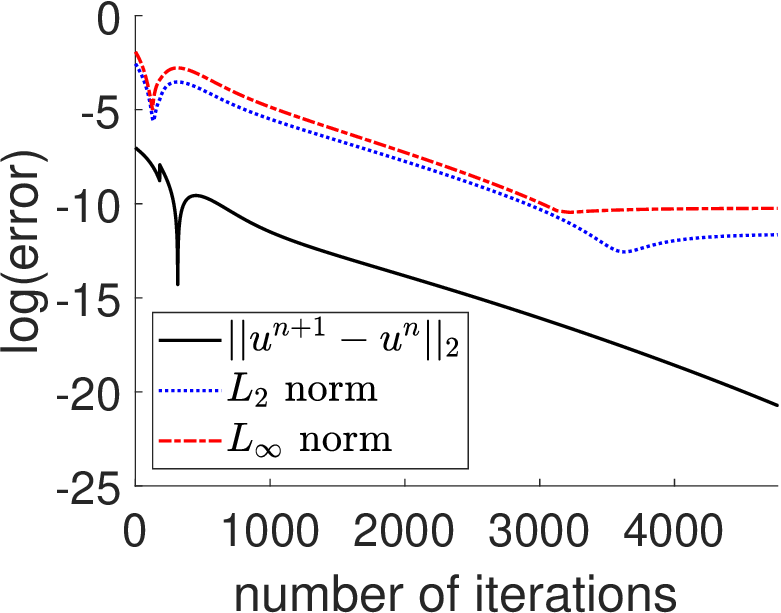}
	\end{tabular}
	
	\caption{(Test problem (\ref{eq.ex1}) with $\alpha=1$. Scheme (\ref{eq.split1Dis.0})-(\ref{eq.split1Dis.2}).) Graphs of the computed solutions and the related convergence history. (a) Regular triangulation of the unit square. (b) Symmetric triangulation of the unit square. (c) Unstructured anisotropic triangulation of the unit square. (d) Unstructured anisotropic triangulation of a half-unit disk. The second-order derivatives are approximated by (\ref{eq.smoothSimp}).}
	\label{fig.ex1.iso}
\end{figure}

\begin{table}[!ht]
	\begin{center}
		(a) \begin{tabular}{|c|c|c||c|c||c|c|}
			\hline
			$h$& Iterations &$\|u^{n+1}-u^n\|$ & $L_2$ norm & rate & $L_{\infty}$ norm& rate \\
			\hline
			1/10 & 193  & 9.88$\times 10^{-10}$  & 6.06$\times 10^{-4}$ &  & 9.86$\times 10^{-4}$ &  \\
			\hline
			1/20 & 606  & 9.81$\times 10^{-10}$  & 1.66$\times 10^{-4}$ & 1.87 & 2.73$\times 10^{-4}$ & 1.85 \\
			\hline
			1/40 & 2064  & 9.96$\times 10^{-10}$  & 4.34$\times 10^{-5}$ & 1.94 & 7.11$\times 10^{-5}$ & 1.94 \\
			\hline
			1/80 & 7577  & 9.99$\times 10^{-10}$  & 1.05$\times 10^{-5}$ & 2.05 & 1.73$\times 10^{-5}$ & 2.04 \\
			\hline
		\end{tabular} \\ \vspace{0.1cm}
		(b) \begin{tabular}{|c|c|c||c|c||c|c|}
			\hline
			$h$& Iterations &$\|u^{n+1}-u^n\|$ & $L_2$ norm & rate & $L_{\infty}$ norm& rate \\
			\hline
			1/10 & 305  & 9.99$\times 10^{-10}$  & 1.37$\times 10^{-3}$ &  & 2.64$\times 10^{-3}$ &  \\
			\hline
			1/20 & 1021  & 9.99$\times 10^{-10}$  & 3.53$\times 10^{-4}$ & 1.96 & 6.88$\times 10^{-4}$ & 1.94 \\
			\hline
			1/40 & 3961  & 9.99$\times 10^{-10}$  & 8.98$\times 10^{-5}$ & 1.97 & 1.75$\times 10^{-4}$ & 1.98 \\
			\hline
			1/80 & 14259  & 9.99$\times 10^{-10}$  & 2.16$\times 10^{-5}$ & 2.06 & 4.24$\times 10^{-5}$ & 2.05 \\
			\hline
		\end{tabular} \\ \vspace{0.1cm}
		(c) \begin{tabular}{|c|c|c||c|c||c|c|}
			\hline
			$h$& Iterations &$\|u^{n+1}-u^n\|$ & $L_2$ norm & rate & $L_{\infty}$ norm& rate \\
			\hline
			1/10 & 180  & 9.38$\times 10^{-10}$  & 5.70$\times 10^{-4}$ &  & 2.04$\times 10^{-3}$ &  \\
			\hline
			1/20 & 591  & 9.80$\times 10^{-10}$  & 1.90$\times 10^{-4}$ & 1.59 & 5.99$\times 10^{-4}$& 1.77 \\
			\hline
			1/40 & 2080  & 9.97$\times 10^{-10}$  & 5.27$\times 10^{-5}$ & 1.85 & 1.57$\times 10^{-4}$ & 1.93 \\
			\hline
			1/80 & 7690  & 9.99$\times 10^{-10}$  & 1.42$\times 10^{-5}$ & 1.89 & 6.22$\times 10^{-5}$ & 1.34 \\
			\hline
		\end{tabular} \\ \vspace{0.1cm}
		(d) \begin{tabular}{|c|c|c||c|c||c|c|}
			\hline
			$h$& Iterations &$\|u^{n+1}-u^n\|$ & $L_2$ norm & rate & $L_{\infty}$ norm& rate \\
			\hline
			1/10 & 111  & 8.49$\times 10^{-10}$  & 6.10$\times 10^{-4}$ &  & 1.20$\times 10^{-3}$ &  \\
			\hline
			1/20 & 374  & 9.97$\times 10^{-10}$  & 1.65$\times 10^{-4}$ & 1.89 & 4.55$\times 10^{-4}$ & 1.40 \\
			\hline
			1/40 & 1221  & 9.93$\times 10^{-10}$  & 3.61$\times 10^{-5}$ & 2.20 & 1.03$\times 10^{-4}$ & 2.14 \\
			\hline
			1/80 & 4765  & 9.97$\times 10^{-10}$  & 8.73$\times 10^{-6}$ & 2.05 & 3.54$\times 10^{-5}$ & 1.54 \\
			\hline
		\end{tabular} \\ \vspace{0.1cm}
	\end{center}
	\caption{(Test problem (\ref{eq.ex1}) with $\alpha=1$. Scheme (\ref{eq.split1Dis.0})-(\ref{eq.split1Dis.2}).) Numbers of iterations necessary for convergence, approximation errors and accuracy orders. (a) Regular triangulation of the unit square. (b) Symmetric triangulation of the unit square. (c) Unstructured anisotropic triangulation of the unit square. (d) Unstructured anisotropic triangulation of the half-unit disk. The second-order derivatives are approximated by (\ref{eq.smoothSimp}).}
	\label{tab.ex1.iso}
\end{table}

\begin{table}[!ht]
	\begin{center}
		(a) \begin{tabular}{|c|c|c||c|c||c|c|}
			\hline
			$h$& Iterations &$\|u^{n+1}-u^n\|$ & $L_2$ norm & rate & $L_{\infty}$ norm& rate \\
			\hline
			1/10 & 266  & 9.50$\times 10^{-8}$  & 1.01$\times 10^{-1}$ &  & 1.22$\times 10^{-1}$ &  \\
			\hline
			1/20 & 512  & 9.87$\times 10^{-8}$  & 4.02$\times 10^{-2}$ & 1.33 & 4.63$\times 10^{-2}$ & 1.40 \\
			\hline
			1/40 & 1432  & 9.99$\times 10^{-8}$  & 1.82$\times 10^{-2}$ & 1.14 & 2.13$\times 10^{-2}$ & 1.12 \\
			\hline
			1/80 & 4529  & 9.99$\times 10^{-8}$  & 8.73$\times 10^{-3}$ & 1.06 & 1.03$\times 10^{-2}$ & 1.05 \\
			\hline
		\end{tabular} \\ \vspace{0.1cm}
		(b) \begin{tabular}{|c|c|c||c|c||c|c|}
			\hline
			$h$& Iterations &$\|u^{n+1}-u^n\|$ & $L_2$ norm & rate & $L_{\infty}$ norm& rate \\
			\hline
			1/10 & 471  & 9.57$\times 10^{-8}$  & 8.24$\times 10^{-2}$ &  & 9.74$\times 10^{-2}$ &  \\
			\hline
			1/20 & 782  & 9.95$\times 10^{-8}$  & 3.46$\times 10^{-2}$ & 1.25 & 3.97$\times 10^{-2}$ & 1.29 \\
			\hline
			1/40 & 2581  & 9.99$\times 10^{-8}$  & 1.60$\times 10^{-2}$ & 1.11 & 1.80$\times 10^{-2}$ & 1.14 \\
			\hline
			1/80 & 7690  & 9.99$\times 10^{-8}$  & 7.78$\times 10^{-3}$ & 1.04 & 8.56$\times 10^{-3}$ & 1.07 \\
			\hline
		\end{tabular} \\ \vspace{0.1cm}
	\end{center}
	\caption{(Test problem (\ref{eq.ex1}) with $\alpha=1$. Scheme (\ref{eq.split1Dis.0})-(\ref{eq.split1Dis.2}).) Numbers of iterations necessary for convergence, approximation errors,  and accuracy orders. (a) Regular triangulation of the unit square. (b) Symmetric triangulation of the unit square. The second-order derivatives are approximated by (\ref{eq.double.6})-(\ref{eq.double.7}).}
	\label{tab.ex1.iso.double}
\end{table}

\begin{table}[!ht]
	\begin{center}
		\begin{tabular}{|c|c|c||c|c||c|c|}
			\hline
			$h$& Iterations &$\|u^{n+1}-u^n\|$ & $L_2$ norm & rate & $L_{\infty}$ norm& rate \\
			\hline
			1/10 & 198  & 9.78$\times 10^{-10}$  & 1.88$\times 10^{-3}$ &  & 2.75$\times 10^{-3}$ &  \\
			\hline
			1/20 & 604  & 9.80$\times 10^{-10}$  & 3.72$\times 10^{-4}$ & 2.34 & 5.91$\times 10^{-4}$& 2.22 \\
			\hline
			1/40 & 2057  & 9.92$\times 10^{-10}$  & 8.86$\times 10^{-5}$ & 2.07 & 1.44$\times 10^{-4}$ & 2.04 \\
			\hline
			1/80 & 7566  & 9.99$\times 10^{-10}$  & 2.14$\times 10^{-5}$ & 2.05 & 3.53$\times 10^{-5}$ & 2.03 \\
			\hline
		\end{tabular} \\ \vspace{0.1cm}
	\end{center}
	\caption{(Test problem (\ref{eq.ex1}) with $\alpha=1$. Scheme (\ref{eq.split2Dis.0})-(\ref{eq.split2Dis.2})) Numbers of iterations necessary for convergence, approximation errors,  and accuracy orders. The second-order derivatives are approximated by (\ref{eq.smoothSimp}).}
	\label{tab.ex1.iso.scheme2}
\end{table}

\begin{table}[!ht]
	\begin{center}
		\begin{tabular}{|c|c|c||c|c||c|c|}
			\hline
			$h$& Iterations &$\|u^{n+1}-u^n\|$ & $L_2$ norm & rate & $L_{\infty}$ norm& rate \\
			\hline
			1/10 & 309  & 9.67$\times 10^{-10}$  & 5.01$\times 10^{-4}$ &  & 8.07$\times 10^{-4}$ &  \\
			\hline
			1/20 & 938  & 9.93$\times 10^{-10}$  & 1.32$\times 10^{-4}$ & 1.92 & 2.12$\times 10^{-4}$ & 1.93 \\
			\hline
			1/40 & 2982  & 9.97$\times 10^{-10}$  & 3.38$\times 10^{-5}$ & 2.01 & 5.39$\times 10^{-5}$ & 1.98 \\
			\hline
			1/80 & 14565  & 9.99$\times 10^{-11}$  & 8.51$\times 10^{-6}$ & 1.99 & 1.36$\times 10^{-5}$ & 1.99 \\
			\hline
		\end{tabular} \\ \vspace{0.1cm}
	\end{center}
	\caption{(Test problem (\ref{eq.ex1}) with $\alpha=2$. Scheme (\ref{eq.split1Dis.0})-(\ref{eq.split1Dis.2}).) Numbers of iterations necessary for convergence, approximation errors,  and accuracy orders on the regular triangulation of the unit square. The second-order derivatives are approximated by (\ref{eq.smoothSimp}).}
	\label{tab.ex1.aniso}
\end{table}

\begin{figure}
	\centering
	\begin{tabular}{ccc}
		(a) & (b) & (c)\\
		\includegraphics[width=0.22\textwidth]{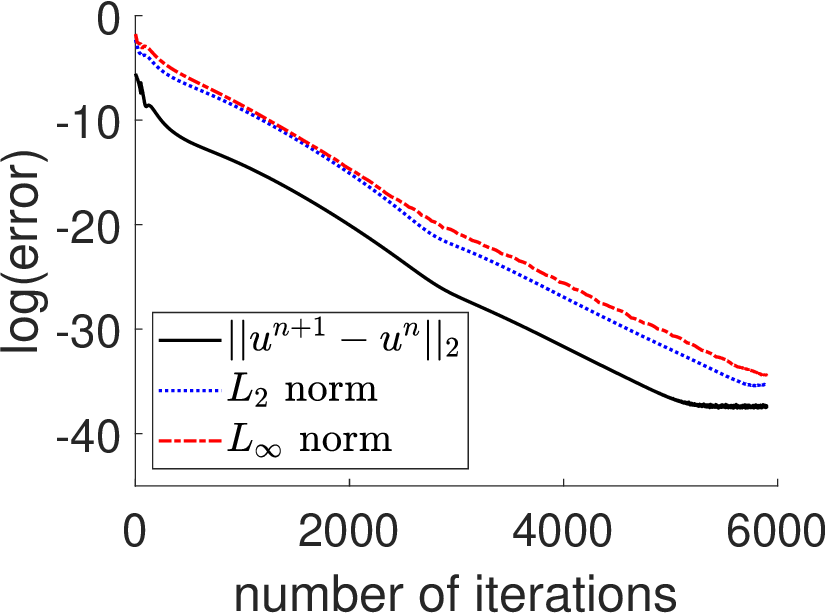} &
		\includegraphics[width=0.22\textwidth]{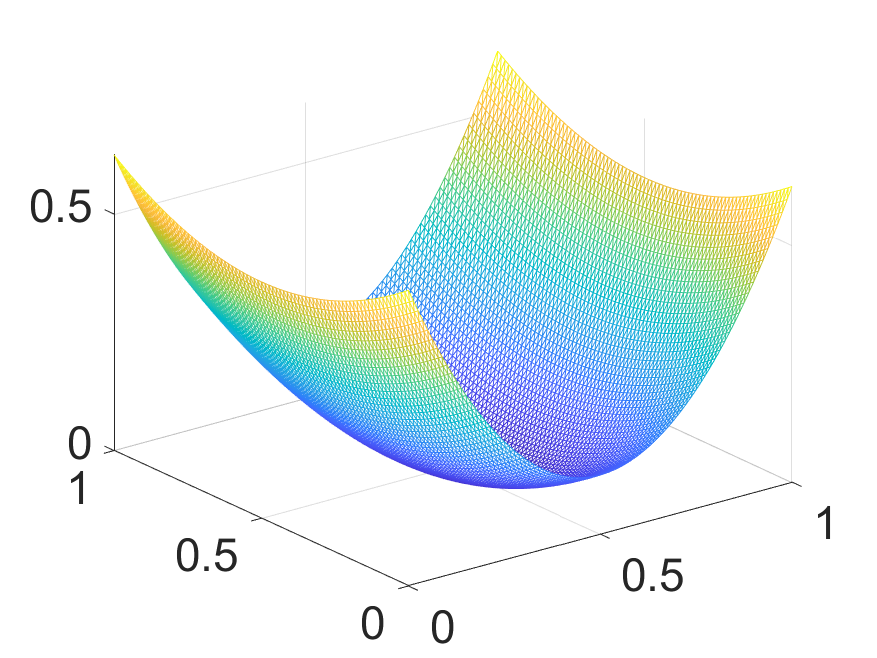}&
		\includegraphics[width=0.22\textwidth]{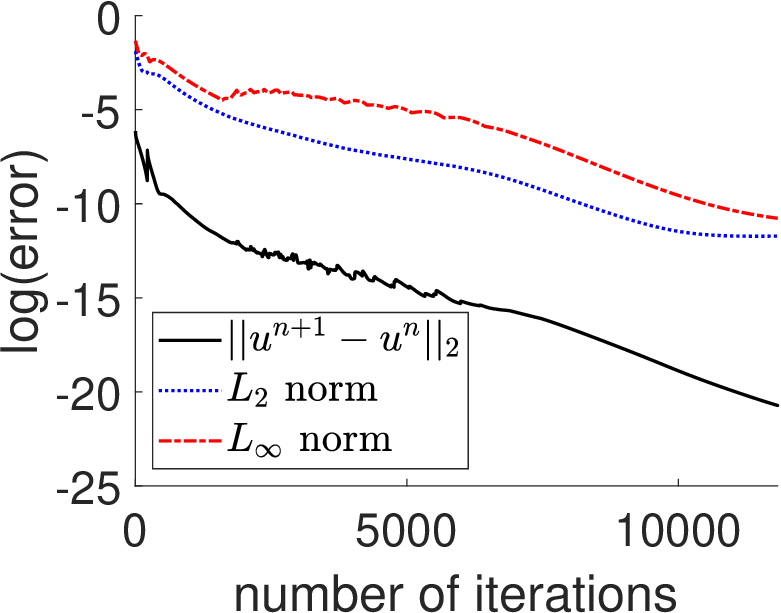}
	\end{tabular}
	
	\caption{(Test problem (\ref{eq.ex1}). Scheme (\ref{eq.split1Dis.0})-(\ref{eq.split1Dis.2}).) (a): With $\alpha=1$, $\varepsilon_1=\varepsilon_2=0$ in both stages, the convergence  history on the regular triangulation of the unit square. The second-order derivatives are approximated by (\ref{eq.smoothSimp}). (b)-(c):The graph of the computed solution and the convergence history on the regular triangulation of the unit square. The second-order derivatives are approximated by (\ref{eq.smoothSimp}).}
	\label{fig.ex1.aniso}
\end{figure}

\subsection{Example 1}
For the first example, we choose the exact solution $u^*$ as a quadratic function,
\begin{equation}
	u^*=\alpha\,(x_1-0.5)^2+ ( x_2-0.5)^2/{\alpha},
	\label{eq.ex1}
\end{equation}
so that  the Gauss curvature
$K=\frac{4}{1+4\alpha\,(x_1-0.5)^2+\frac{4}{\alpha}(x_2-0.5)^2}$ 
and the boundary condition $g=u^*|_{\partial\Omega}$, where $\alpha$ is a positive constant. 

Since the solution of this example is smooth, we use (\ref{eq.gradApprox}) to approximate the first-order derivatives. In the first test, we choose $\alpha=1$ so that $u^*$ represents a family of con-centric circles which vary isotropically.

With the second-order derivatives approximated by (\ref{eq.smoothSimp}) and scheme (\ref{eq.split1Dis.0})-(\ref{eq.split1Dis.2}), the graphs and convergence histories of numerical solutions on different meshes are shown in Figure \ref{fig.ex1.iso}. The numbers of iteration and accuracy orders are shown in Table \ref{tab.ex1.iso}, where the accuracy orders in the $L_2$ and $L_{\infty}$ norms are in general larger than $1.5$. Stopping criterion $\|u^{n+1}-u^n\|<10^{-9}$ is used.

With the second-order derivatives approximated by (\ref{eq.double.6}) and (\ref{eq.double.7}) and scheme (\ref{eq.split1Dis.0})-(\ref{eq.split1Dis.2}), we can use a less demanding stopping criterion. Here we use $\|u^{n+1}-u^n\|<10^{-7}$. The numbers of iteration, the errors of approximation, and the rates of convergence on regular and symmetric meshes of the unit square are shown in Table \ref{tab.ex1.iso.double}, which demonstrate that, in general, our algorithm with approximation (\ref{eq.double.6}) and (\ref{eq.double.7}) is first-order accurate, and in comparison with the results based on the approximation (\ref{eq.smoothSimp}), the errors based on the approximation  (\ref{eq.double.6}) and (\ref{eq.double.7}) are larger and the convergence rates are smaller.

For comparison, we also show the results by scheme (\ref{eq.split2Dis.0})-(\ref{eq.split2Dis.2}) with the second-order derivatives approximated by (\ref{eq.smoothSimp}). Since the boundary condition is compatible, we use a large $\kappa=500$. The number of iteration and accuracy orders are shown in Table \ref{tab.ex1.iso.scheme2}. Its efficiency and accuracy are similar to that of scheme (\ref{eq.split1Dis.0})-(\ref{eq.split2Dis.2}). If $\kappa$ goes to infinity, scheme (\ref{eq.split2Dis.0})-(\ref{eq.split2Dis.2}) has an additional stabilization (diffusion) term which provides larger error but extra robustness, the same as what is observed by comparing Table \ref{tab.ex1.iso}(a) and Table \ref{tab.ex1.iso.scheme2}.

Since the exact solution is a quadratic function, its second-order derivatives are constants so that the zero Neumann boundary condition on these derivatives is exact. With $\varepsilon_1=\varepsilon_2=0$ and $h=1/40$, the convergence history of scheme (\ref{eq.split1Dis.0})-(\ref{eq.split1Dis.2}) is shown in Figure \ref{fig.ex1.aniso}(a). We can see that although approximation (\ref{eq.smoothSimp}) is a kind of variational crime, the error decreases to machine precision.

In the second test, we choose $\alpha=2$ so that $u^*$ represents a family of con-centric ellipses which vary anisotropically. We apply our algorithm to this problem on the unit square with regular meshes. The number of iterations necessary to satisfy the stopping criterion and the corresponding approximation error accuracy is shown in Table \ref{tab.ex1.aniso}. The graph of the computed solution and the related convergence history are shown in Figure \ref{fig.ex1.aniso}(b)-(c).

\begin{center}
	\begin{table}[!hbp]
		\centering
		(a)
		\begin{tabular}{|c|c|c|c|c|c|c|}
			\hline
			$h$&   Iteration & $\|u^{n+1}-u^n\|_2$ & $L_2$ error &  rate & $L_{\infty}$ error & rate \\
			\hline
			1/16 & 177& 9.57$\times 10^{-7}$ & 9.79$\times 10^{-2}$ &  & 1.69$\times 10^{-1}$ &  \\
			\hline
			1/32 & 791& 9.98$\times 10^{-7}$ & 5.61$\times 10^{-2}$ & 0.80 & 1.19$\times 10^{-1}$ & 0.51 \\
			\hline
			1/64 & 3360& 9.97$\times 10^{-7}$ & 3.12$\times 10^{-2}$ & 0.85 & 8.39$\times 10^{-2}$ & 0.50 \\
			\hline
			1/128 & 17273& 9.99$\times 10^{-7}$ & 1.55$\times 10^{-2}$ & 1.01 & 5.86$\times 10^{-2}$ & 0.52 \\
			\hline
		\end{tabular}\\\vspace{0.1cm}
		(b)
		\begin{tabular}{|c|c|c|c|c|c|c|c|c|}
			\hline
			$h$&   Iteration & $\|u^{n+1}-u^n\|_2$ & $L_2$ error &  rate & $L_{\infty}$ error & rate \\
			\hline
			1/16 & 236& 9.72$\times 10^{-7}$ & 2.86$\times 10^{-2}$ &  & 7.40$\times 10^{-2}$ &  \\
			\hline
			1/32 & 1179& 9.98$\times 10^{-7}$ & 1.13$\times 10^{-2}$ & 1.34 & 4.53$\times 10^{-2}$ & 0.71 \\
			\hline
			1/64 & 5261& 9.95$\times 10^{-7}$ & 7.33$\times 10^{-3}$ & 0.62 & 4.12$\times 10^{-2}$ & 0.14 \\
			\hline
		\end{tabular}\\\vspace{0.1cm}
		(c)
		\begin{tabular}{|c|c|c|}
			\hline
			$h$& $L_{\infty}$ error & rate\\
			\hline
			1/16 & 1.61$\times 10^{-1}$ & \\
			\hline
			1/32 &  1.28$\times 10^{-1}$& 0.33\\
			\hline
			1/64 &  1.09$\times 10^{-1}$& 0.23\\
			\hline
			1/128 & 8.80$\times 10^{-2}$ & 0.31\\
			\hline
		\end{tabular}\\\vspace{0.1cm}
		\caption{(Test problem (\ref{eq.ex3.1}). Scheme (\ref{eq.split1Dis.0})-(\ref{eq.split1Dis.2}).) Numbers of iterations, approximation errors, and accuracy orders with the second-order derivatives approximated by (a) (\ref{eq.smoothSimp}) and (b) (\ref{eq.double.6})-(\ref{eq.double.7}). (c) shows the $L_{\infty}$ errors and accuracy orders from \cite{hamfeldt2018convergent}.}
		\label{tab.ex3.1}
	\end{table}
\end{center}

\begin{figure}
	\centering
	\begin{tabular}{cccc}
		(a) & (b) & (c) & (d)\\
		\includegraphics[width=0.21\textwidth]{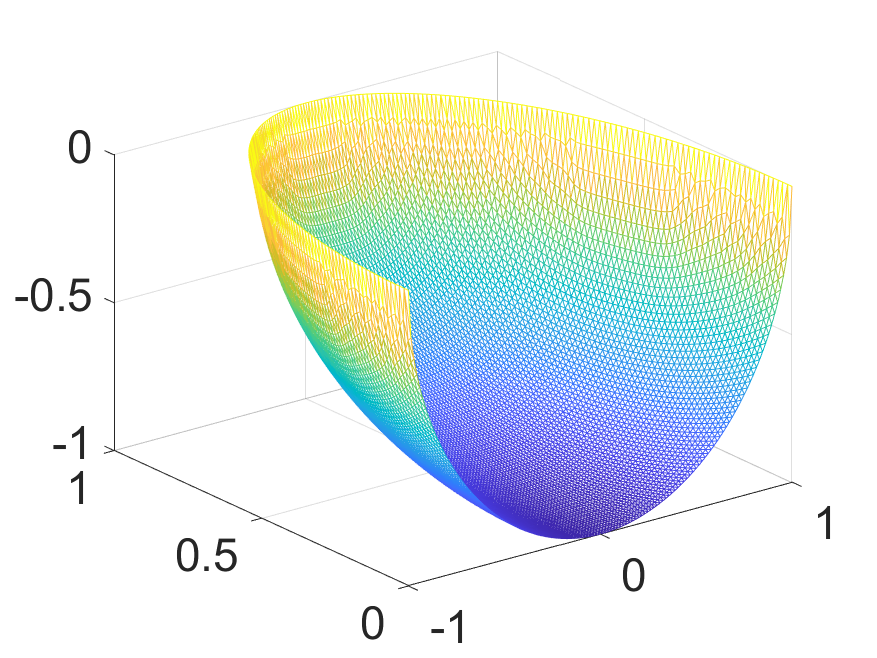}&
		\includegraphics[width=0.21\textwidth]{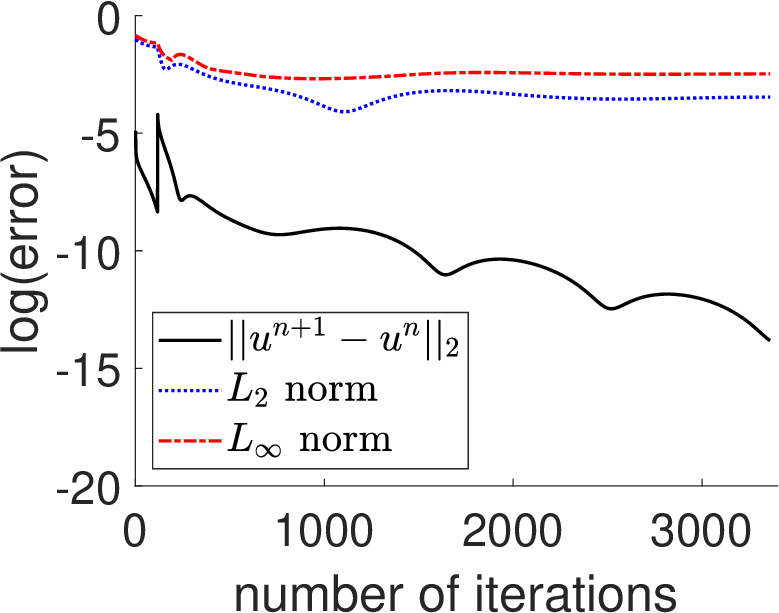}&
		\includegraphics[width=0.21\textwidth]{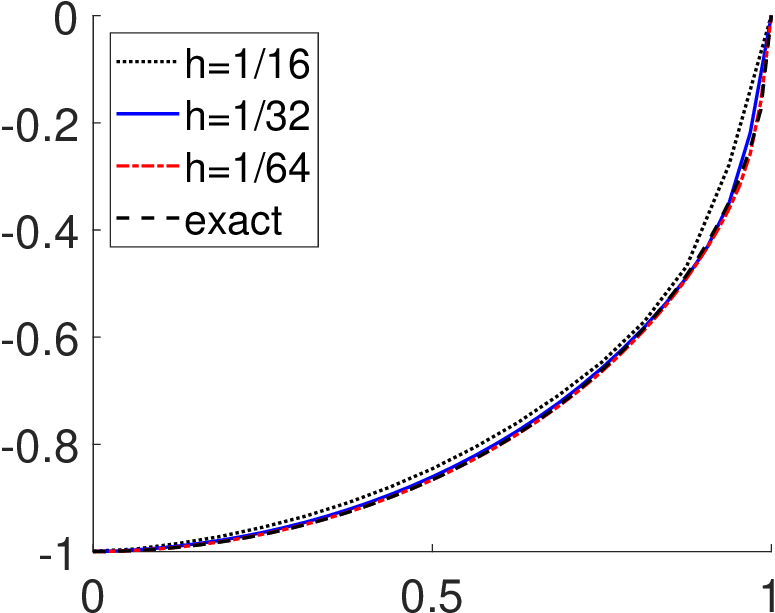} &
		\includegraphics[width=0.21\textwidth]{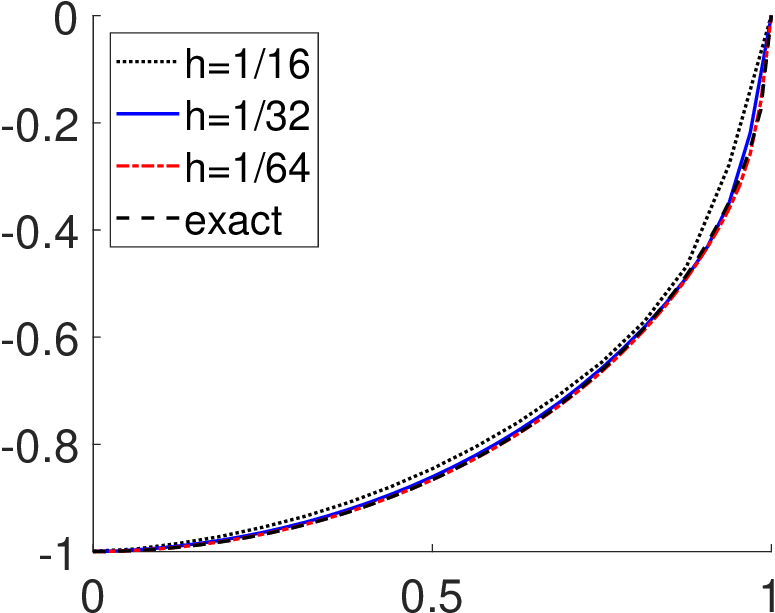}
	\end{tabular}
	
	\caption{(Test problem (\ref{eq.ex3.1}). Scheme (\ref{eq.split1Dis.0})-(\ref{eq.split1Dis.2})) (a) The graph of the solution with $h=1/64$, where the second-order derivatives are approximated by (\ref{eq.smoothSimp}). (b) The convergence history of (a). (c)-(d): Graphs of the restrictions of the numerical solutions to the line $x_1=0$ with different $h$'s. (c) the second-order derivatives approximated by (\ref{eq.smoothSimp}) and (d) the second-order derivatives approximated by (\ref{eq.double.6})-(\ref{eq.double.7}).}
	\label{fig.ex3.vc.sol}
\end{figure}

\subsection{Example 2}\label{sec.sphere1}
In the second example, we consider a problem with the exact solution
\begin{equation}
	u=-\sqrt{1-x_1^2-x_2^2},
	\label{eq.ex3.1}
\end{equation}
which is a part of the unit sphere, and the corresponding Gauss curvature is constant: $K=1$. The computational domain is chosen to be half of the unit disk, $\Omega=\{(x_1,x_2)\,|\, x_1\geq0,\, x_1^2+x_2^2\leq 1\}$. Accordingly, the boundary condition is given as
\begin{equation}
	g=
	\begin{cases}
		0, & x_1>0, \\
		-\sqrt{1-x_2^2}, & x_1=0.
	\end{cases}
\end{equation}

This problem is interesting since the gradient of the exact solution along the boundary where $x_1>0$ is infinite (a more challenging problem is solved in Section \ref{sec.sphere12}). This problem is also solved in \cite{hamfeldt2018convergent}.  Since the first-order derivatives are infinite along a part of the boundary, we have to use the regularized approximation (\ref{eq.firstapprox0}) for the first-order derivatives; otherwise our solution will blow up. We use $\varepsilon=h$, $\varepsilon_1=h^2$, $\Delta t=h^2$, and the stopping criterion $\|u^{n+1}-u^n\|_2<10^{-6}$. Figure \ref{fig.ex3.vc.sol}(a)-(b) shows the graph of the numerical solution for $h=1/64$ and the related convergence history with second-order derivatives approximated by (\ref{eq.smoothSimp}). The cross sections of the numerical solutions along the boundary $x_1=0$ with second-order derivatives approximated by (\ref{eq.smoothSimp}) or (\ref{eq.double.6})-(\ref{eq.double.7}) are shown in Figure \ref{fig.ex3.vc.sol}(c)-(d), and the convergence of numerical solutions using both approximations is  clearly observed.

To further quantify both approximations of second-order derivatives, we show the numbers of iterations, the $L_2$- and $L_{\infty}$- errors, and their corresponding rates of convergence in Tables \ref{tab.ex3.1}(a) and (b). In Table \ref{tab.ex3.1}, we can see that both approximations of the second-order derivatives behave reasonably well. Although the algorithm equipped with approximation (\ref{eq.double.6})-(\ref{eq.double.7}) produces smaller errors than the one equipped with (\ref{eq.smoothSimp}), the algorithm with approximation (\ref{eq.smoothSimp}) is more stable as its convergence rate is uniformly about $0.5$. As a comparison, we also list the $L_{\infty}$-errors and related convergence rates from \cite{hamfeldt2018convergent} in Table \ref{tab.ex3.1}(c). When the mesh is fine enough, our algorithm equipped with either  approximations produces smaller $L_{\infty}$-errors than that of \cite{hamfeldt2018convergent}.

\begin{table}
	(a)\\
	\begin{tabular}{|c|c|c|c|c|c|c|}
		\hline
		\scriptsize{$h$}& \scriptsize{Iter.} & \scriptsize{$\|u^{n+1}-u^n\|_2$} &\scriptsize{$\|\p^n-\D^2u^n\|_2$} & \scriptsize{$\frac{\|\p^n-\D^2u^n\|_2}{\|\p^n\|_2}$}& \scriptsize$\min$ & \scriptsize{$\|\p^n-\D^2u^n\|_2$ in.}\\
		\hline
		1/20 & 177& 9.62$\times 10^{-7}$ & 4.44$\times 10^{-2}$ & 2.35$\times 10^{-2}$&-0.1192 & $3.28\times 10^{-3}$\\
		\hline
		1/40 & 672& 9.72$\times 10^{-7}$ & 1.80$\times 10^{-1}$& 7.26$\times 10^{-2}$ & -0.1263& $5.98\times 10^{-3}$\\
		\hline
		1/80 & 2149& 9.98$\times 10^{-7}$ & 4.94$\times 10^{-1}$ & 1.67$\times 10^{-1}$& -0.1305& $9.41\times 10^{-3}$\\
		\hline
	\end{tabular}\\\vspace{0.1cm}
	
	(b) \\
	\begin{tabular}{|c|c|c|c|c|c|c|}
		\hline
		\scriptsize{$h$}& \scriptsize{Iter.} & \scriptsize$\|u^{n+1}-u^n\|_2$ & \scriptsize{$\|\p^n-\D^2u^n\|_2$} & \scriptsize{$\frac{\|\p^n-\D^2u^n\|_2}{\|\p^n\|_2}$}& \scriptsize$\min$ & \scriptsize{$\|\p^n-\D^2u^n\|_2$ in.}\\
		\hline
		1/20 & 246& 9.88$\times 10^{-7}$ & 4.60$\times 10^{-3}$ & 2.43$\times 10^{-3}$&-0.1345 & $7.47\times 10^{-4}$\\
		\hline
		1/40 & 695& 9.88$\times 10^{-7}$ & 7.53$\times 10^{-2}$& 2.43$\times 10^{-2}$ & -0.1359& $4.77\times 10^{-3}$\\
		\hline
		1/80 & 2468& 9.99$\times 10^{-7}$ & 3.52$\times 10^{-1}$ & 6.82$\times 10^{-2}$& -0.1376& $6.54\times 10^{-3}$\\
		\hline
	\end{tabular}\\\vspace{0.1cm}
	\caption{(Test problem (\ref{eq.ex5}). Scheme (\ref{eq.split1Dis.0})-(\ref{eq.split1Dis.2}).) Numbers of iterations, iteration errors and minimum values. The second-order derivatives are approximated by (a) (\ref{eq.smoothSimp}) and (b) (\ref{eq.double.6})-(\ref{eq.double.7}).}
	\label{tab.ex5}
\end{table}

\begin{figure}
	\centering
	\begin{tabular}{cccc}
		(a) & (b) & (c) & (d)\\
		\includegraphics[width=0.21\textwidth]{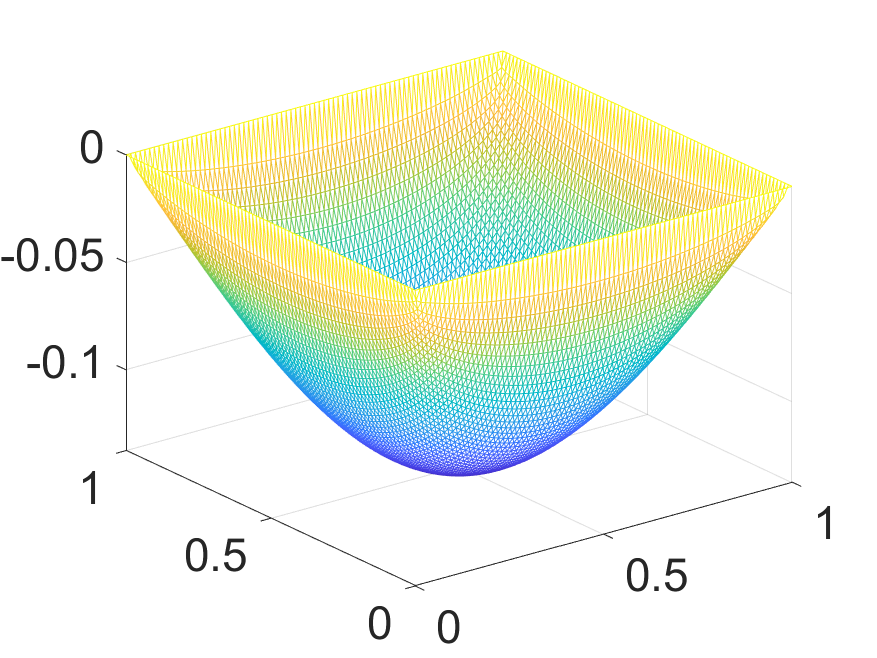}&
		\includegraphics[width=0.21\textwidth]{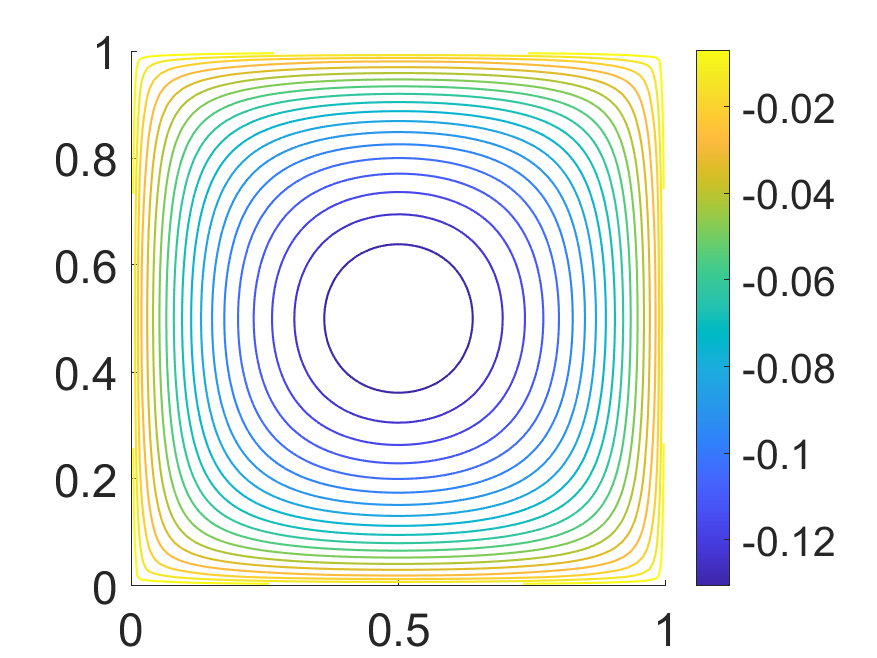}&
		\includegraphics[width=0.21\textwidth]{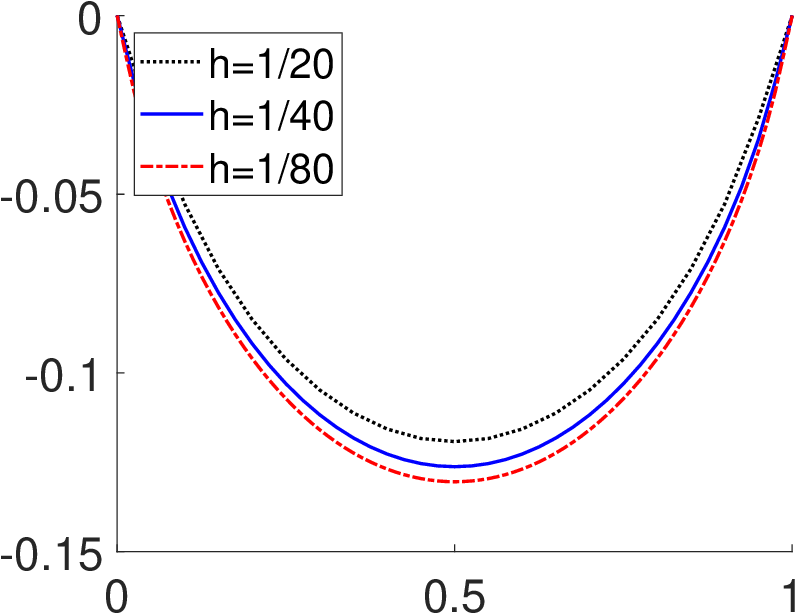}&
		\includegraphics[width=0.21\textwidth]{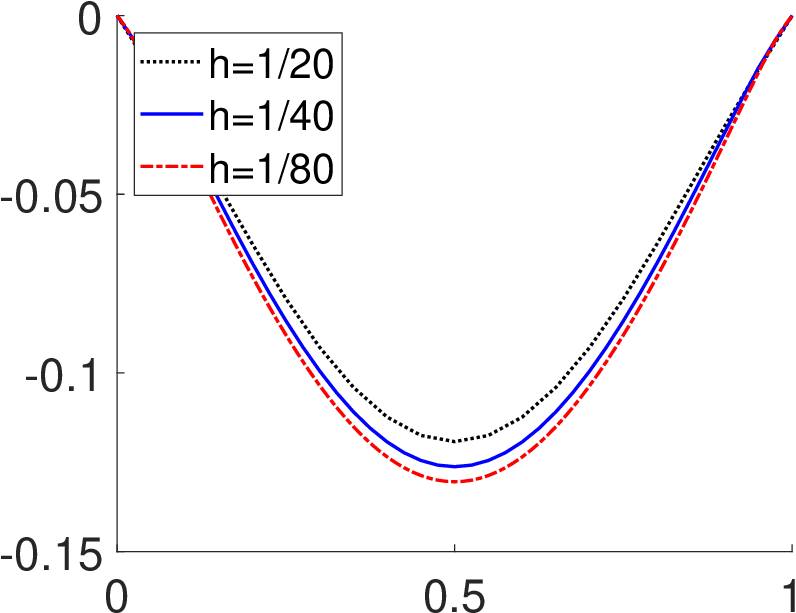}\\
		& (e) & (f) &\\
		&\includegraphics[width=0.21\textwidth]{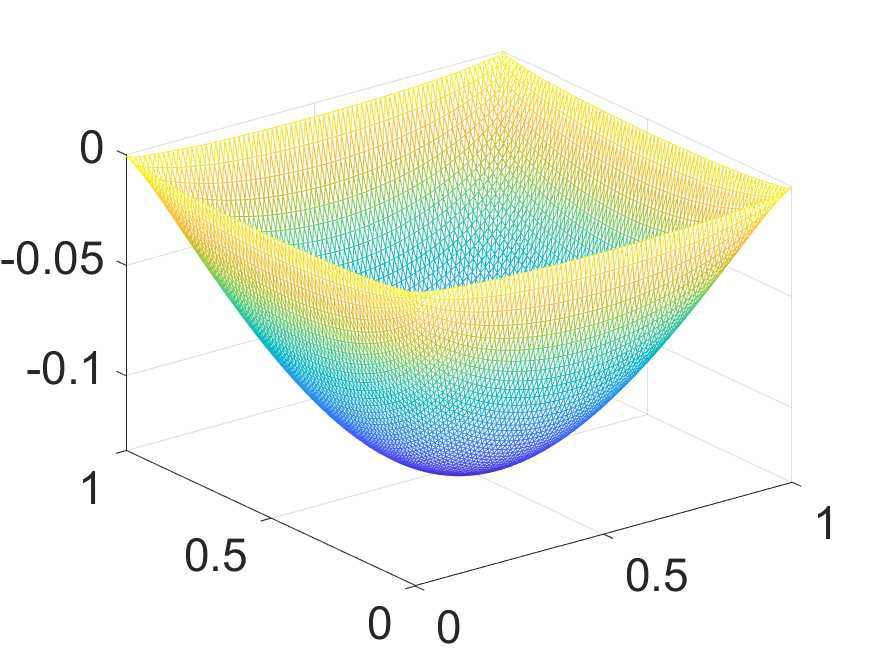}&
		\includegraphics[width=0.21\textwidth]{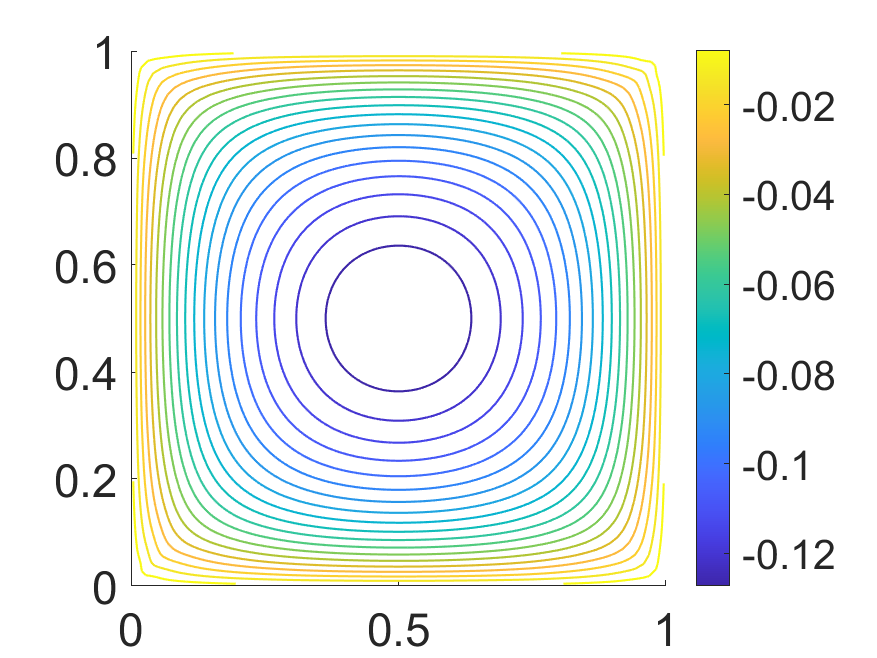}
	\end{tabular}
	
	\caption{(Test problem (\ref{eq.ex5}).) (a)-(d):  Scheme (\ref{eq.split1Dis.0})-(\ref{eq.split1Dis.2}). (a)-(b) Graphs and contour of the numerical solution with $h=1/80$. (c)-(d) Graphs of the restrictions of numerical solutions along $x_1=1/2$ and $x_1=x_2$ with $h=1/20$, $1/40$, and  $1/80$. 
		(e)-(f): Scheme (\ref{eq.split2Dis.0})-(\ref{eq.split2Dis.2}).) (e) Graphs of the numerical solution with $h=1/80$. (f) Contour of the solution displayed in (a). The second-order derivatives are approximated by (\ref{eq.smoothSimp}).}
	\label{fig.ex5}
\end{figure}

\subsection{Example 3}\label{sec.feng.square}
We consider a problem with no classical solution. The curvature is a constant in $\Omega$:
\begin{equation}
	K=1/2 \mbox{ in } \Omega,
	\label{eq.ex5}
\end{equation}
where $\Omega=[0,1]^2$. We use the boundary condition $g=0 \mbox{ on } \partial \Omega$. This problem has no classical solution since $\det(\D^2u)$ vanishes on $\partial\Omega$. In other words, this problem is incompatible. In our experiment, we first use scheme (\ref{eq.split1Dis.0})-(\ref{eq.split1Dis.2}) with $\varepsilon=\varepsilon_1=h^2$ and  $\Delta t=2h^2$. The second-order derivatives are approximated by (\ref{eq.smoothSimp}). The number of iteration, convergence errors, and minimum values are shown in Table \ref{tab.ex5}. The graphs and contour of the numerical solution with $h=1/80$ are shown in Figure \ref{fig.ex5} (a) and (b). The comparisons of the restriction of the numerical solution with different $h$ along $x_1=1/2$ and $x_1=x_2$ are shown in Figure \ref{fig.ex5} (c) and (d). Our solution is smooth and almost convex, except the region near the corners of the domain.

Then we use scheme (\ref{eq.split2Dis.0})-(\ref{eq.split1Dis.2}) with $\varepsilon=\varepsilon_1=h^2$ and $\Delta t=8h^2$ to solve it. With $h=1/80$, the graph and contour of the computed solution are shown in Figure \ref{fig.ex5}. We can see the boundary value of the computed solution is no longer constant. At the middle segment on each edges, its values is away from 0 to be compatible with its interior value. The same problem is solved on an ellipse domain in Section \ref{sec.feng.ellipse}.

\subsection{Example 4}
\begin{figure}
	\centering
	\begin{tabular}{cccc}
		(a) & (b) & (c) & (d)\\
		\includegraphics[width=0.21\textwidth]{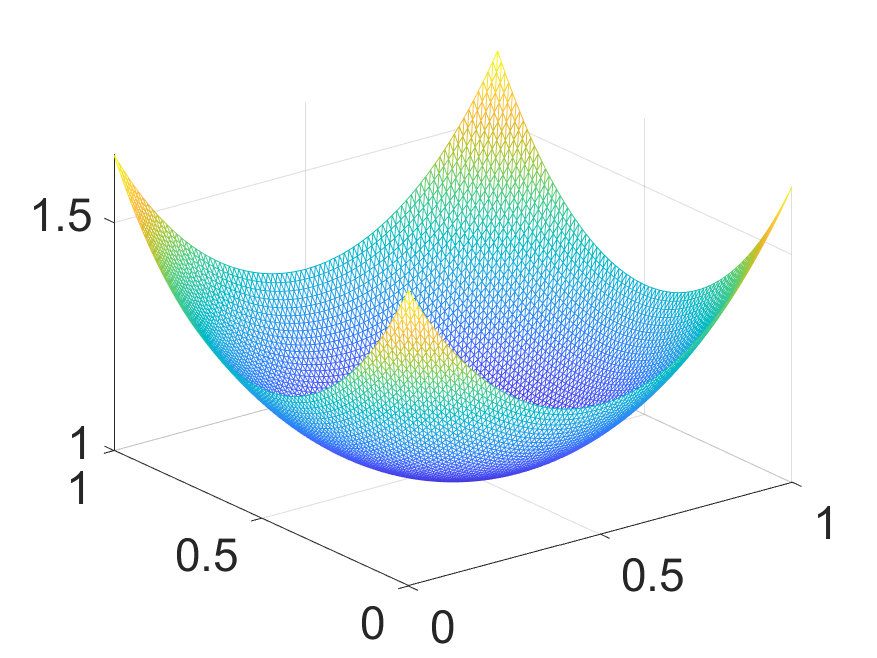} & 
		\includegraphics[width=0.21\textwidth]{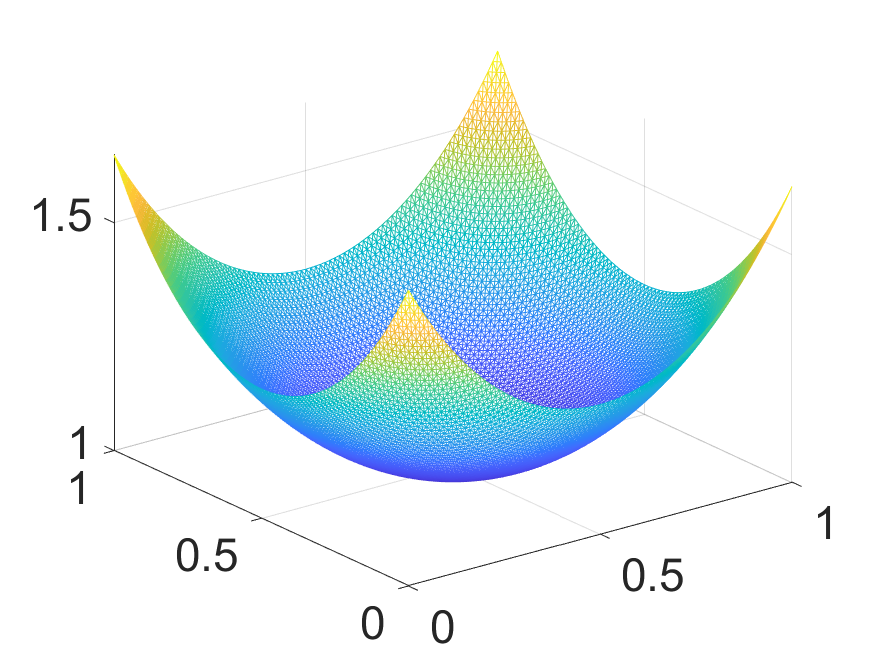}&
		\includegraphics[width=0.21\textwidth]{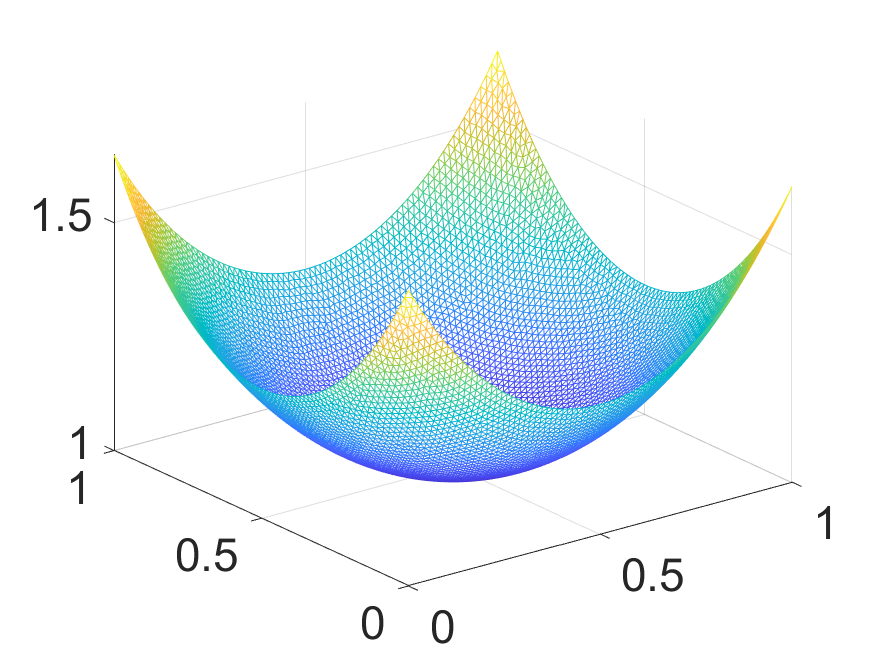}&
		\includegraphics[width=0.21\textwidth]{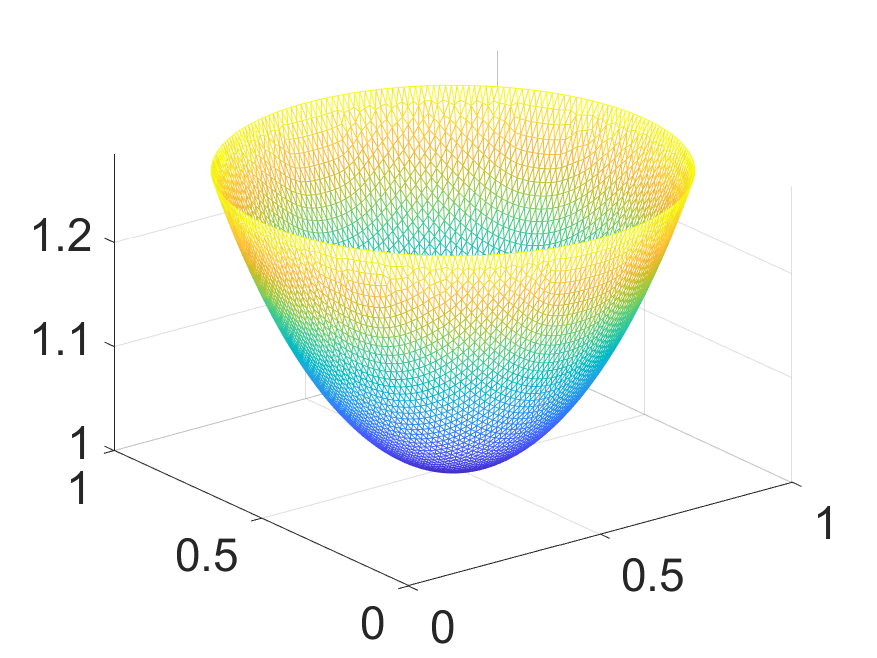}\\
		\includegraphics[width=0.21\textwidth]{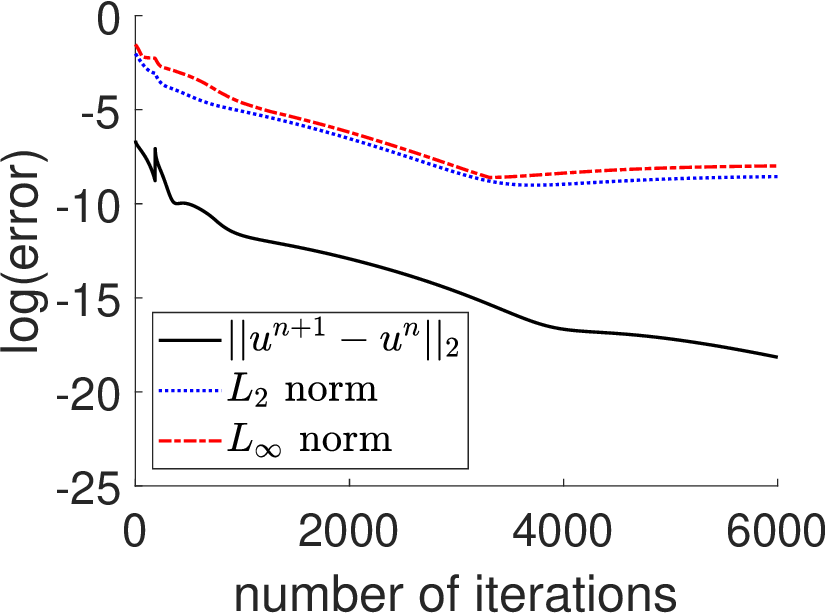}&
		\includegraphics[width=0.21\textwidth]{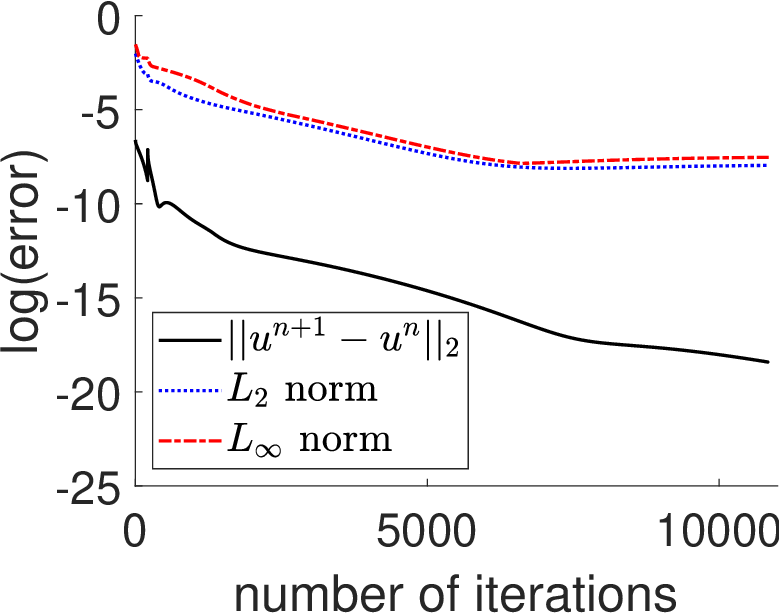}&
		\includegraphics[width=0.21\textwidth]{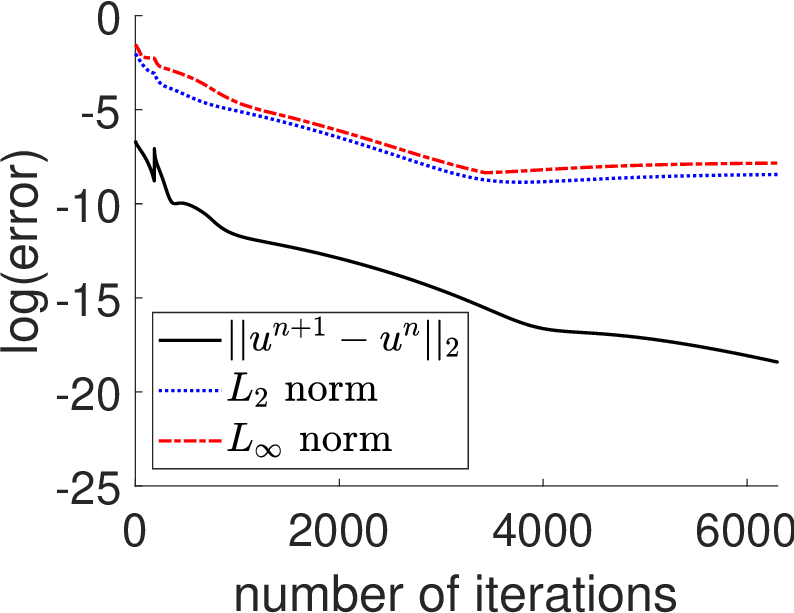}&
		\includegraphics[width=0.21\textwidth]{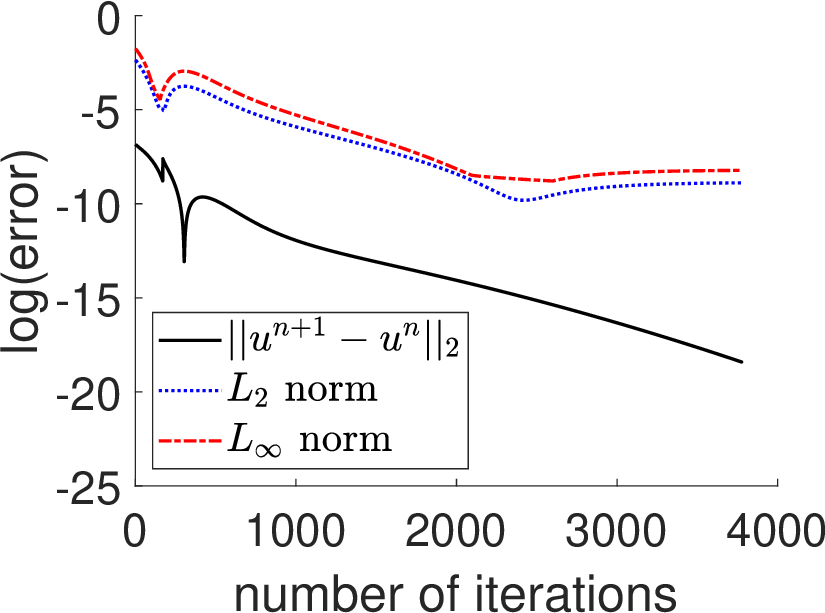}
	\end{tabular}

	\caption{(Test problem (\ref{eq.ex2}). Scheme (\ref{eq.split1Dis.0})-(\ref{eq.split1Dis.2}).) Graphs of the computed solutions and the related convergence history. (a) Regular triangulation of the unit square. (b) Symmetric triangulation of the unit square. (c) Unstructured anisotropic triangulation of the unit square. (d) Unstructured anisotropic triangulation of the half-unit disk. The second order derivatives are approximated by (\ref{eq.smoothSimp}).}
	\label{fig.ex2}
\end{figure}

In this example, we choose the exact solution $u^*$ as an exponential function:
\begin{equation}
	u^*=e^{r^2},
	\label{eq.ex2}
\end{equation}
so that $g=u^*|_{\Omega}$ and the Gauss curvature
$$
K=\frac{4e^{2r^2}(1+2r^2)}{(1+4r^2e^{2r^2})^2},
$$
where $r^2=(x_1-1/2)^2+(x_2-1/2)^2$  in (\ref{eq.ex2}). For this particular example, we choose the stopping criterion $\|u^{n+1}-u^n\|_2<10^{-8}$. The graphs of numerical solutions and convergence histories are shown in Figure \ref{fig.ex2}. In Table \ref{tab.ex2},  we show the numbers of iterations, the $L_2$- and $L_{\infty}$- errors, and the convergence orders for different meshes with various mesh sizes, and as we can see, the convergence orders for the four different meshes are close to $2$ if the meshes are fine enough.

\begin{table}[!ht]
	\begin{center}
		(a) \begin{tabular}{|c|c|c||c|c||c|c|}
			\hline
			$h$& Iterations &$\|u^{n+1}-u^n\|$ & $L_2$ norm & rate & $L_{\infty}$ norm& rate \\
			\hline
			1/10 & 180  & 9.25$\times 10^{-9}$  & 8.04$\times 10^{-3}$ &  & 1.52$\times 10^{-2}$ &  \\
			\hline
			1/20 & 550  & 9.99$\times 10^{-9}$  & 2.78$\times 10^{-3}$ & 1.53 & 4.98$\times 10^{-3}$ & 1.61 \\
			\hline
			1/40 & 1825  & 9.96$\times 10^{-9}$  & 7.73$\times 10^{-4}$ & 1.85 & 1.36$\times 10^{-3}$ & 1.87 \\
			\hline
			1/80 & 6222  & 9.99$\times 10^{-9}$  & 1.95$\times 10^{-4}$ & 1.99 & 3.42$\times 10^{-4}$ & 1.99 \\
			\hline
		\end{tabular} \\ \vspace{0.1cm}
		(b) \begin{tabular}{|c|c|c||c|c||c|c|}
			\hline
			$h$& Iterations &$\|u^{n+1}-u^n\|$ & $L_2$ norm & rate & $L_{\infty}$ norm& rate \\
			\hline
			1/10 & 326  & 9.86$\times 10^{-9}$  & 1.11$\times 10^{-2}$ & & 1.95$\times 10^{-2}$ &  \\
			\hline
			1/20 & 949  & 9.97$\times 10^{-9}$  & 4.00$\times 10^{-3}$ & 1.47 & 6.64$\times 10^{-3}$ & 1.55 \\
			\hline
			1/40 & 3346  & 9.99$\times 10^{-9}$  & 1.25$\times 10^{-3}$ & 1.68 & 1.99$\times 10^{-3}$ & 1.74 \\
			\hline
			1/80 & 10845  & 9.99$\times 10^{-9}$  & 3.49$\times 10^{-4}$ & 1.84 & 5.35$\times 10^{-4}$ & 1.90 \\
			\hline
		\end{tabular} \\ \vspace{0.1cm}
		(c) \begin{tabular}{|c|c|c||c|c||c|c|}
			\hline
			$h$& Iterations &$\|u^{n+1}-u^n\|$ & $L_2$ norm & rate & $L_{\infty}$ norm& rate \\
			\hline
			1/10 & 166  & 9.65$\times 10^{-9}$  & 8.25$\times 10^{-3}$ & & 1.63$\times 10^{-2}$ &  \\
			\hline
			1/20 & 537  & 9.82$\times 10^{-9}$  & 3.04$\times 10^{-3}$ & 1.44 & 5.32$\times 10^{-3}$ & 1.62 \\
			\hline
			1/40 & 1837 & 9.94$\times 10^{-9}$  & 8.42$\times 10^{-4}$ & 1.85 & 1.53$\times 10^{-3}$ & 1.80 \\
			\hline
			1/80 & 6303 & 9.99$\times 10^{-9}$  & 2.14$\times 10^{-4}$ & 1.98 & 3.93$\times 10^{-4}$ & 1.96 \\
			\hline
		\end{tabular} \\ \vspace{0.1cm}
		(d) \begin{tabular}{|c|c|c||c|c||c|c|}
			\hline
			$h$& Iterations &$\|u^{n+1}-u^n\|$ & $L_2$ norm & rate & $L_{\infty}$ norm& rate \\
			\hline
			1/10 & 95  & 8.34$\times 10^{-9}$  & 5.28$\times 10^{-3}$ &  & 1.14$\times 10^{-2}$ &  \\
			\hline
			1/20 & 310 & 9.91$\times 10^{-9}$  & 1.79$\times 10^{-3}$ & 1.56 & 3.66$\times 10^{-3}$ & 1.64 \\
			\hline
			1/40 & 1005  & 9.93$\times 10^{-9}$  & 4.99$\times 10^{-4}$ & 1.84 & 9.94$\times 10^{-4}$ & 1.88 \\
			\hline
			1/80 & 3779  & 9.99$\times 10^{-9}$  & 1.38$\times 10^{-4}$ & 1.85 & 2.67$\times 10^{-4}$ & 1.90 \\
			\hline
		\end{tabular} \\ \vspace{0.1cm}
	\end{center}
	\caption{(Test problem (\ref{eq.ex2}). Scheme (\ref{eq.split1Dis.0})-(\ref{eq.split1Dis.2}).) Numbers of iterations necessary for convergence, approximation errors, and accuracy orders. (a) Regular triangulation of the unit square. (b) Symmetric triangulation of the unit square. (c) Unstructured anisotropic triangulation of the unit square. (d) Unstructured anisotropic triangulation of the half-unit disk. The second-order derivatives are approximated by (\ref{eq.smoothSimp}).}
	\label{tab.ex2}
\end{table}

\begin{center}
	\begin{table}
		(a)\\
		\begin{tabular}{|c|c|c|c|c|c|c|c|}
			\hline
			$h$&   Iter & $\|u^{n+1}-u^n\|_2$ & $L_2$ norm &  rate & $L_{\infty}$ norm & rate& min\\
			\hline
			1/8 & 366& 9.79$\times 10^{-7}$ & 9.64$\times 10^{-2}$ &  & 1.26$\times 10^{-1}$ & &-1.0744\\
			\hline
			1/16 & 1323& 9.85$\times 10^{-7}$ & 3.94$\times 10^{-2}$ & 1.29 & 4.38$\times 10^{-2}$ & 1.52 & -1.0426\\
			\hline
			1/32 & 3189& 9.98$\times 10^{-7}$ & 1.59$\times 10^{-2}$ & 1.31 & 1.79$\times 10^{-2}$ & 1.29 & -1.0312\\
			\hline
		\end{tabular}\\\vspace{0.1cm}
		(b)\\
		\begin{tabular}{|c|c|c|c|c|c|c|c|}
			\hline
			$h$&   Iter & $\|u^{n+1}-u^n\|_2$ & $L_2$ norm &  rate & $L_{\infty}$ norm & rate& min\\
			\hline
			1/8 & 171& 9.39$\times 10^{-7}$ & 4.11$\times 10^{-2}$ &  & 4.20$\times 10^{-2}$ & &-0.4524\\
			\hline
			1/16 & 845& 9.93$\times 10^{-7}$ & 2.42$\times 10^{-2}$ & 0.76 & 2.12$\times 10^{-2}$ & 0.99 & -0.4408\\
			\hline
			1/32 & 3636& 9.99$\times 10^{-7}$ & 1.14$\times 10^{-2}$ & 1.09 & 1.01$\times 10^{-2}$ & 1.07 & -0.4311\\
			\hline
		\end{tabular}\\\vspace{0.1cm}
		(c)\\
		\begin{tabular}{|c|c|c|c|c|c|c|c|}
			\hline
			$h$&   Iter & $\|u^{n+1}-u^n\|_2$ & $L_2$ norm &  rate & $L_{\infty}$ norm & rate& min\\
			\hline
			1/8 & 120& 9.20$\times 10^{-7}$ & 5.09$\times 10^{-2}$ &  & 3.89$\times 10^{-2}$ & &-0.2281\\
			\hline
			1/16 & 634& 9.94$\times 10^{-7}$ & 2.91$\times 10^{-2}$ & 0.81 & 2.20$\times 10^{-2}$ & 0.82 & -0.2138\\
			\hline
			1/32 & 2829& 9.99$\times 10^{-7}$ & 1.55$\times 10^{-2}$ & 0.91 & 1.18$\times 10^{-2}$ & 0.90 & -0.2039\\
			\hline
		\end{tabular}\\\vspace{0.1cm}
		(d)\\
		\begin{tabular}{|c|c|c|c|c|c|c|c|}
			\hline
			$h$&   Iter & $\|u^{n+1}-u^n\|_2$ & $L_2$ norm &  rate & $L_{\infty}$ norm & rate& min\\
			\hline
			1/8 & 154& 9.41$\times 10^{-7}$ & 7.38$\times 10^{-2}$ &  & 5.54$\times 10^{-2}$ & &-0.1452\\
			\hline
			1/16 & 510& 9.96$\times 10^{-7}$ & 4.38$\times 10^{-2}$ & 0.75 & 3.18$\times 10^{-2}$ & 0.80 & -0.1258\\
			\hline
			1/32 & 2445& 9.97$\times 10^{-7}$ & 2.51$\times 10^{-2}$ & 0.80 & 1.82$\times 10^{-2}$ & 0.81 & -0.1125\\
			\hline
		\end{tabular}\\\vspace{0.1cm}
		\caption{(Test problem (\ref{eq.ex4}). Scheme (\ref{eq.split1Dis.0})-(\ref{eq.split1Dis.2}).) Numbers of iterations, approximation errors, convergence rates, and minimum values with (a) $R=1/2+0.1$, (b) $R=1$, (c) $R=2$, (d) $R=4$. The exact minimum values for (a)-(d) are $-1.0221$, $-0.4226$, $-0.1926$, and $-0.0944$, respectively. The second-order derivatives are approximated by (\ref{eq.smoothSimp}).}
		\label{tab.ex4}
	\end{table}
\end{center}

\begin{figure}
	\centering
	\begin{tabular}{cccc}
		& (a) & (b) & \\
		&\includegraphics[width=0.21\textwidth]{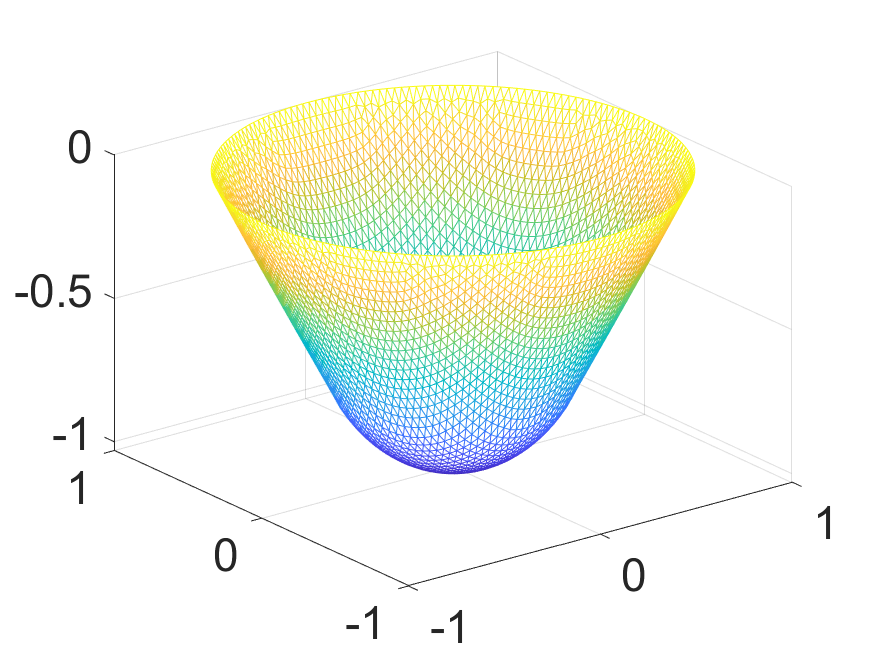}&
		\includegraphics[width=0.21\textwidth]{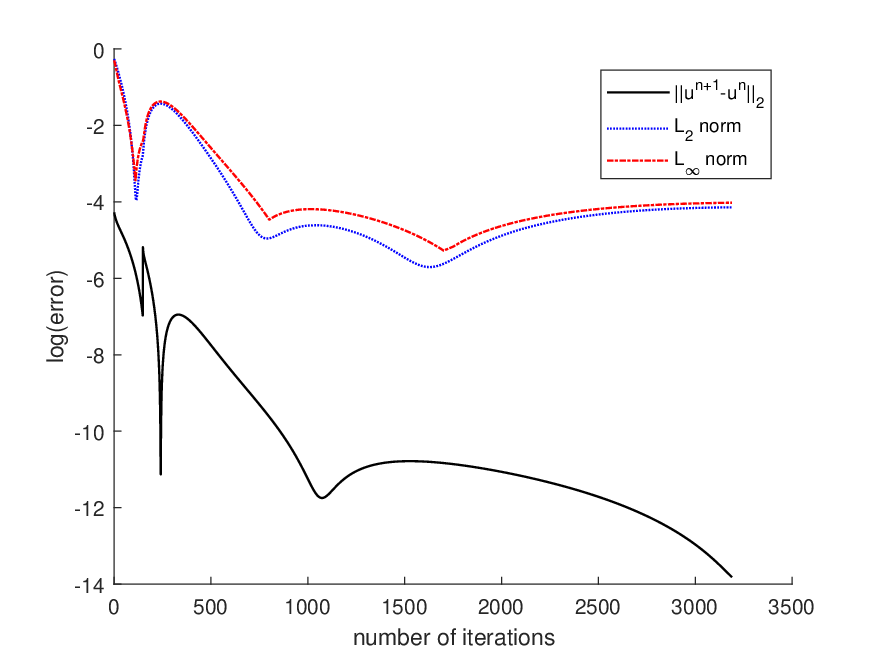}&\\
		(c) & (d) & (e) & (f)\\
		\includegraphics[width=0.21\textwidth]{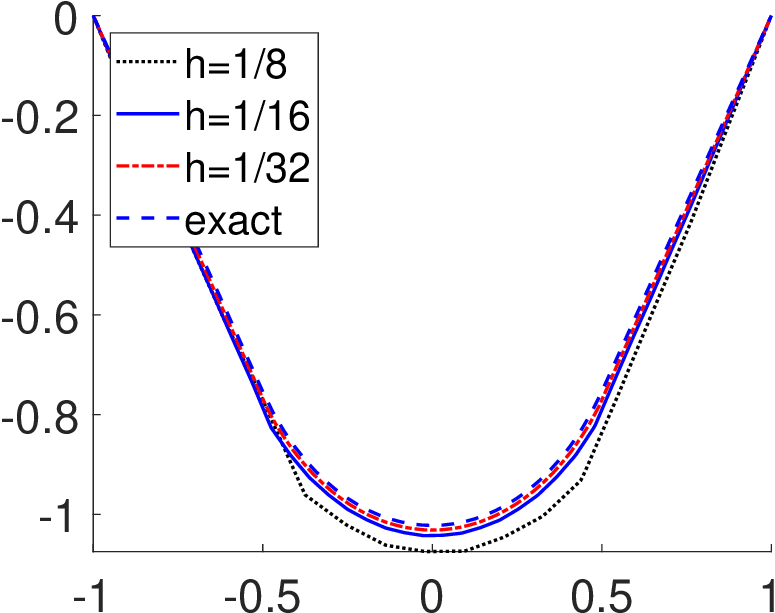}&
		\includegraphics[width=0.21\textwidth]{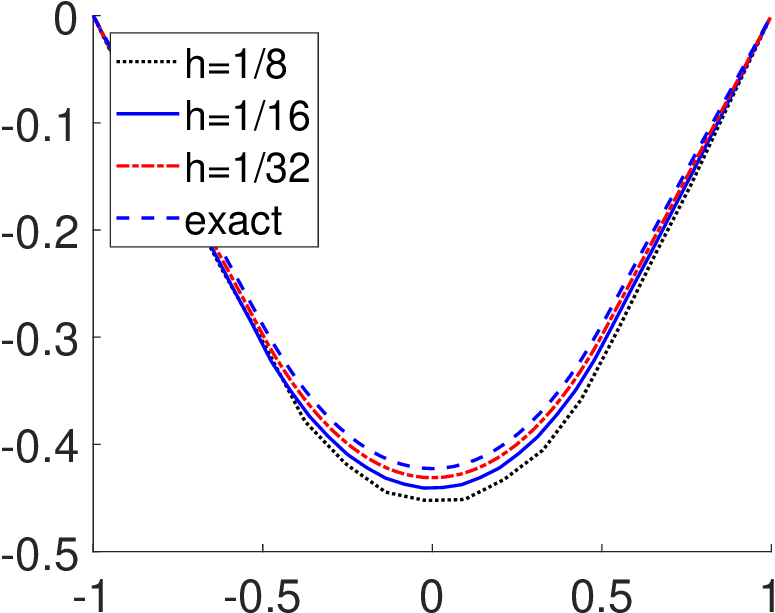}&
		\includegraphics[width=0.21\textwidth]{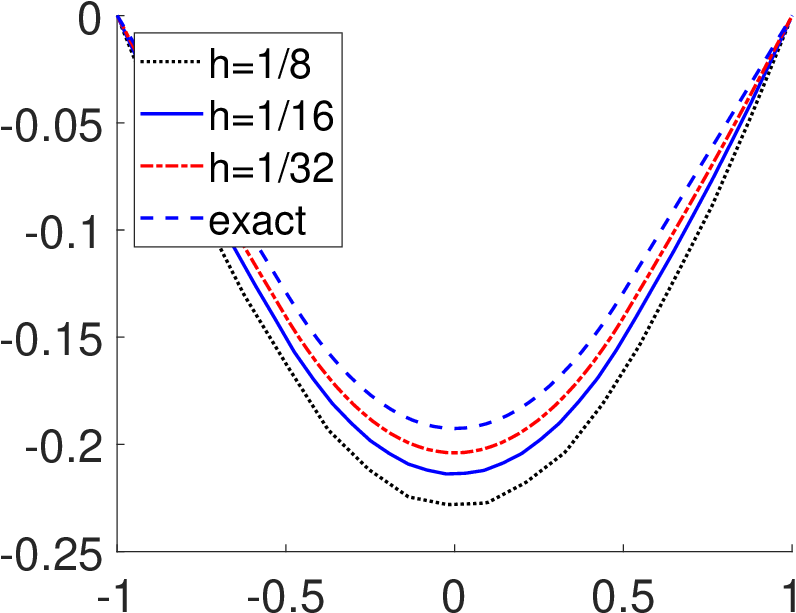}&
		\includegraphics[width=0.21\textwidth]{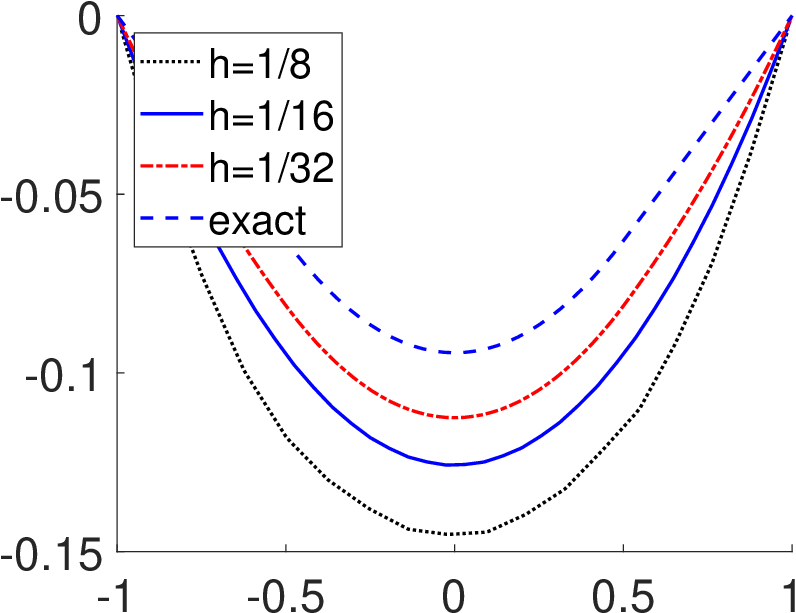}
	\end{tabular}
	
	\caption{(Test problem (\ref{eq.ex4}). Scheme (\ref{eq.split1Dis.0})-(\ref{eq.split1Dis.2}).) (a)-(b): The graph of the numerical solution and the convergence history with $h=1/32$ and $R=1/2+0.1$. The second-order derivatives are approximated by (\ref{eq.smoothSimp}). (c)-(f): Graphs of the restrictions of numerical solutions along $x_1=0$ with $h=1/8$, $1/16$, and  $1/32$ and (c) $R=1/2+0.1$, (d) $R=1$, (e) $R=2$, and (f) $R=4$. The exact minimum values for (a)-(d) are $-1.0221$, $-0.4226$, $-0.1926$, and $-0.0944$, respectively. The second-order derivatives are approximated by (\ref{eq.smoothSimp}). }
	\label{fig.ex4.s}
\end{figure}

\subsection{Example 5}
In this example, the Gauss curvature is given as
\begin{equation}
	K(x)=
	\begin{cases}
		R^{-2} &\mbox{ if } x\in \omega,\\
		\eta & \mbox{ if } x\in\Omega\slash\bar{\omega},
	\end{cases}
	\label{eq.ex4}
\end{equation}
where $\eta$ is close to 0, $R>1/2$, $\Omega=\{x=(x_1,x_2), x_1^2+x_2^2<1\}$, and $\omega=\{x=(x_1,x_2),x_1^2+x_2^2<1/4\}$. The Dirichlet boundary condition is specified as $g=0$ on $\partial \Omega$.

The above Minkowski problem has the exact solution
\begin{equation}
	u^*=\begin{cases}
		\frac{2R^2-1}{\sqrt{4R^2-1}}-\sqrt{R^2-r^2} & \mbox{ if } 0\leq r\leq1/2,\\
		\frac{1}{\sqrt{4R^2-1}}(r-1) & \mbox{ if }1/2\leq r\leq1, 
	\end{cases}
	\label{eq.ex4.exact}
\end{equation}
where, for each $R$, the minimum value of $u^*$ is $\frac{2R^2-1}{\sqrt{4R^2-1}}-R.$ By the form of the exact solution  (\ref{eq.ex4.exact}), we can see that $u^*$ extends to $-\infty$ as $R\rightarrow 1/2$.

We test our algorithm by taking $R=1/2+0.1$, $1$, $2$ and $4$, respectively. In our experiments, we set $\eta=h$ in the Minkowski problem and choose  $\varepsilon=\varepsilon_1=h^2$ and $\Delta t=4h^2$ in our algorithm. We further approximate the first-order derivatives by (\ref{eq.gradApprox}) and the second-order derivatives by (\ref{eq.smoothSimp}).

Figure \ref{fig.ex4.s} shows the graph of the numerical solution and the convergence history with $R=1/2+0.1$ and $h=1/32$.  Figure \ref{fig.ex4.s} shows the restrictions of the numerical solutions along the line $x_1=0$ with different mesh sizes and  different $R$'s. Since we use $\eta=h$, the first-order convergence rate is expected. Table \ref{tab.ex4} shows the number of iterations necessary for convergence to steady state, minimum values, $L_2$- and $L_{\infty}$- errors, and convergence rates for different $R$'s. As expected, the convergence rate is close to first order, while the rate is even higher than first order  for $R=1/2+0.1$.

\subsection{Example 6}\label{sec.sphere12}

We consider the problem which is a half sphere defined on
$$\Omega=\{(x_1,x_2)|\; \sqrt{(x_1-1/2)^2+(x_2-1/2)^2}\leq1/2\}$$
so that the exact solution is
\begin{equation}
	u^*=-\sqrt{\frac{1}{4}-\left(x_1-\frac{1}{2}\right)^2-\left(x_2-\frac{1}{2}\right)^2}.
	\label{eq.ex3.2}
\end{equation}
Accordingly, we have the Gauss curvature $K=4$ and the Dirichlet boundary condition $g=0$ on $\partial\Omega$.

\begin{figure}
	\centering
	\begin{tabular}{cccc}
		(a) & (b) & (c) & (d)\\
		\includegraphics[width=0.21\textwidth]{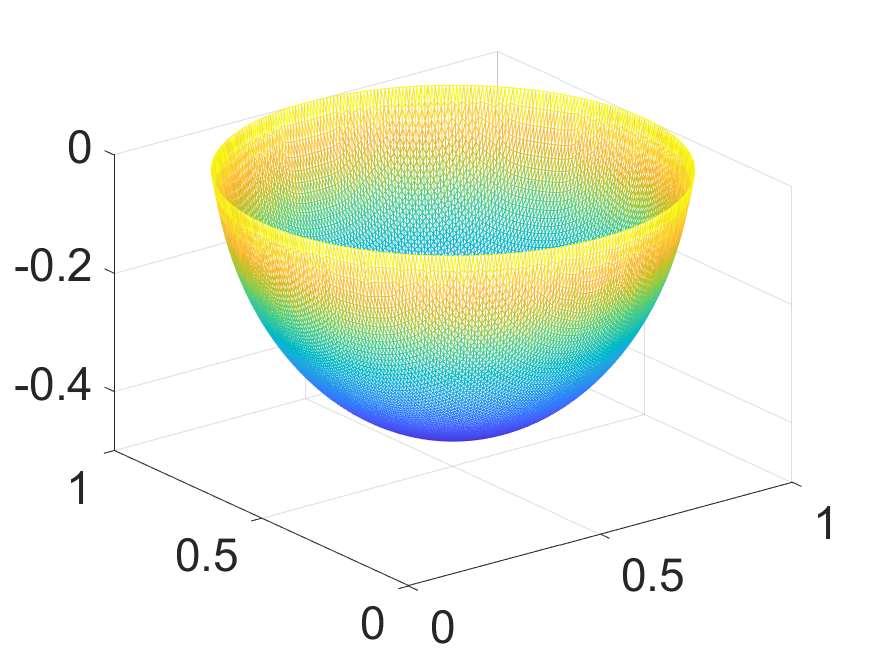}&
		\includegraphics[width=0.21\textwidth]{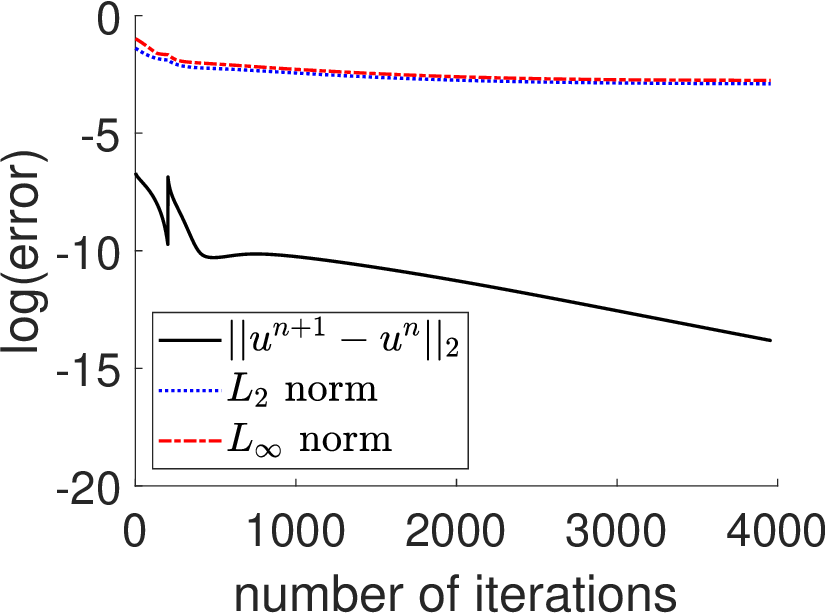}&
		\includegraphics[width=0.21\textwidth]{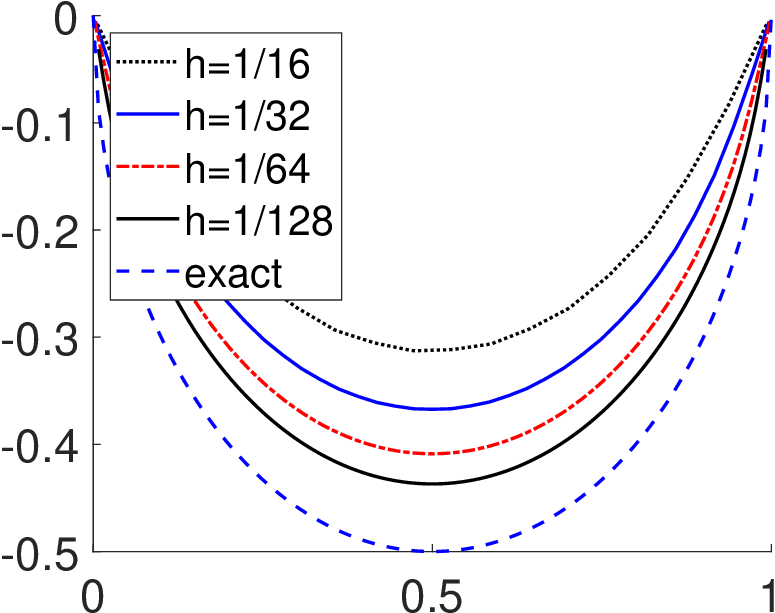}&
		\includegraphics[width=0.21\textwidth]{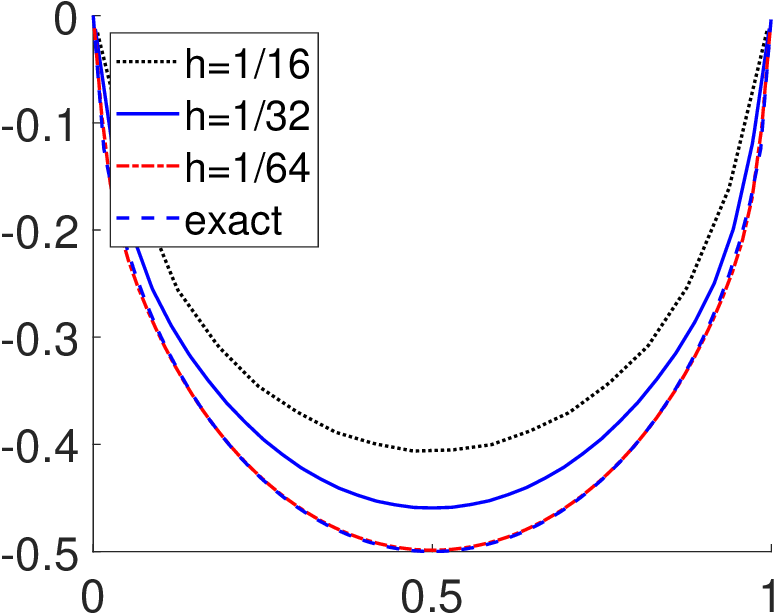}
	\end{tabular}
	
	\caption{(Test problem (\ref{eq.ex3.2}). Scheme (\ref{eq.split1Dis.0})-(\ref{eq.split1Dis.2}).) (a)-(b): Graphs of the approximated solution and the convergence history with $h=1/128$, where  the second-order derivatives are approximated by (\ref{eq.smoothSimp}). (c)-(d): Graphs of the restrictions of the approximated solutions along $x_1=1/2$ with the
		second-order derivatives approximated by (c) (\ref{eq.smoothSimp}) and (d) (\ref{eq.double.6})-(\ref{eq.double.7}).}
	\label{fig.ex3.2s}
\end{figure}

This problem is more challenging than Example (\ref{eq.ex3.1}) in Section \ref{sec.sphere1} because its first-order derivatives along the entire boundary $\partial\Omega$ are infinite. Still, we need to use (\ref{eq.firstapprox0}) to approximate the first-order derivatives. Figure \ref{fig.ex3.2s} shows the graph of the numerical solution and the convergence history, where $h=1/128$ and second-order derivatives are approximated by (\ref{eq.smoothSimp}). Figure \ref{fig.ex3.2s} shows the restrictions of the numerical solutions along the line $x_1=1/2$, where  different mesh sizes are used and second-order derivatives are approximated by (\ref{eq.smoothSimp})
or (\ref{eq.double.6})-(\ref{eq.double.7}). Table \ref{tab.ex3.2} shows the minimum values of the numerical solutions, $L_2$- and $L_{\infty}$-errors, and convergence rates on different meshes, where the minimum value of the exact solution is $-0.5$. Similar to our observation for problem (\ref{eq.ex3.1}), our algorithm with approximation (\ref{eq.double.6})-(\ref{eq.double.7}) produces  smaller errors but the performance of our algorithm with approximation (\ref{eq.smoothSimp}) is more stable. On the fine mesh with $h=1/128$, our algorithm with approximation (\ref{eq.double.6})-(\ref{eq.double.7}) does not converge. Nevertheless, our algorithm with approximation (\ref{eq.smoothSimp}) has convergence rate $0.5$ in both the $L_2$- and $L_{\infty}$- norms.

\begin{center}
	\begin{table}
		(a)\\
		\begin{tabular}{|c|c|c|c|c|c|c|c|}
			\hline
			$h$&   Iter & $\|u^{n+1}-u^n\|_2$ & $L_2$ norm &  rate & $L_{\infty}$ norm & rate& min\\
			\hline
			1/16 & 42& 8.76$\times 10^{-7}$ & 1.47$\times 10^{-1}$ &  & 1.86$\times 10^{-1}$ & &-0.313\\
			\hline
			1/32 & 204& 9.94$\times 10^{-7}$ & 1.10$\times 10^{-1}$ & 0.42 & 1.32$\times 10^{-1}$ & 0.49 & -0.367\\
			\hline
			1/64 & 856& 9.96$\times 10^{-7}$ & 7.85$\times 10^{-2}$ & 0.49 & 9.12$\times 10^{-2}$ & 0.53 & -0.409\\
			\hline
			1/128 & 3956& 9.99$\times 10^{-7}$ & 5.52$\times 10^{-2}$ & 0.51 & 6.36$\times 10^{-2}$ & 0.52 & -0.437\\
			\hline
		\end{tabular}\\\vspace{0.1cm}
		(b)\\
		\begin{tabular}{|c|c|c|c|c|c|c|c|}
			\hline
			$h$&   Iter & $\|u^{n+1}-u^n\|_2$ & $L_2$ norm &  rate & $L_{\infty}$ norm & rate& min\\
			\hline
			1/16 & 59& 8.24$\times 10^{-7}$ & 6.77$\times 10^{-2}$ &  & 9.30$\times 10^{-2}$ & &-0.406\\
			\hline
			1/32 & 307& 9.68$\times 10^{-7}$ & 3.07$\times 10^{-2}$ & 1.14 & 4.80$\times 10^{-2}$ & 0.95 & -0.459\\
			\hline
			1/64 & 1679& 9.99$\times 10^{-7}$ & 4.08$\times 10^{-3}$ & 2.91 & 2.07$\times 10^{-2}$ & 1.21 & -0.499\\
			\hline
		\end{tabular}\\\vspace{0.1cm}
		\caption{(Test problem (\ref{eq.ex3.2}). Scheme (\ref{eq.split1Dis.0})-(\ref{eq.split1Dis.2}).) Numbers of iterations necessary for convergence, approximation errors, and rates of convergence with the second-order derivatives approximated by (a) (\ref{eq.smoothSimp}) and (b) (\ref{eq.double.6})-(\ref{eq.double.7}).}
		\label{tab.ex3.2}
	\end{table}
\end{center}

\subsection{Example 7}
\label{sec.feng.ellipse}
\begin{table}
	\begin{tabular}{|c|c|c|c|c|c|}
		\hline
		$h$&   Iteration & $\|u^{n+1}-u^n\|_2$ &  $\|\p^n-\D^2u^n\|_2$ & $\frac{\|\p^n-\D^2u^n\|_2}{\|\p^n\|_2}$& min\\
		\hline
		1/20 & 324& 9.98$\times 10^{-7}$ & 2.08$\times 10^{-2}$ & 8.81$\times 10^{-3}$&-0.1912\\
		\hline
		1/40 & 1151& 9.93$\times 10^{-7}$ & 4.77$\times 10^{-2}$& 1.99$\times 10^{-2}$ & -0.1926\\
		\hline
		1/80 & 3907& 9.99$\times 10^{-7}$ & 9.61$\times 10^{-2}$ & 4.08$\times 10^{-2}$& -0.1933\\
		\hline
	\end{tabular}\\\vspace{0.1cm}
	\caption{(Test problem (\ref{eq.ex5ellipse}). Scheme (\ref{eq.split1Dis.0})-(\ref{eq.split1Dis.2}).) Numbers of iterations, iteration errors and minimum values. The second-order derivatives are approximated by (\ref{eq.smoothSimp}).}
	\label{tab.ex5ellipse}
\end{table}
\begin{figure}
	\centering
	\begin{tabular}{cccc}
		(a) & (b) & (c) & (d)\\
		\includegraphics[width=0.21\textwidth]{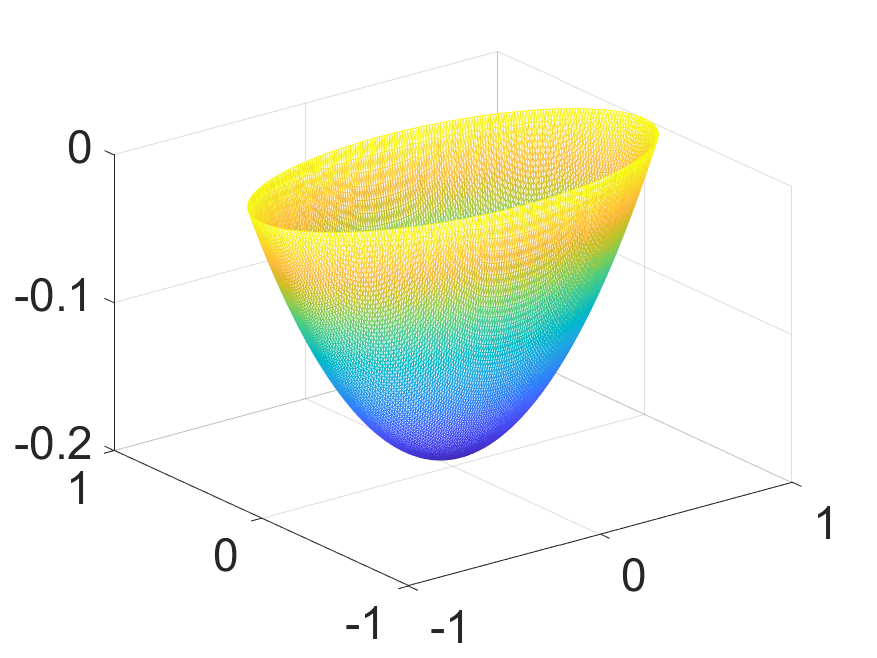}&
		\includegraphics[width=0.21\textwidth]{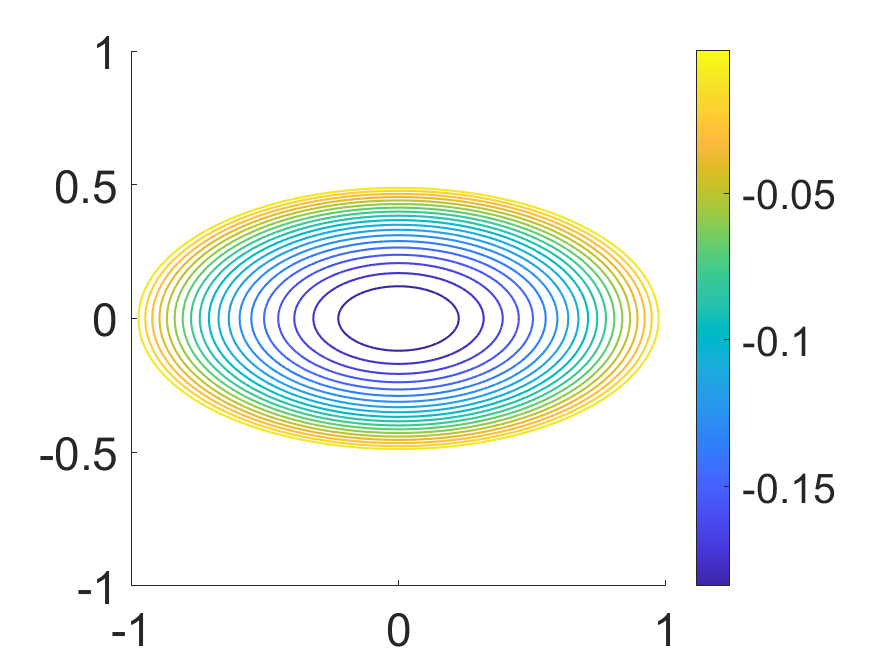}&
		\includegraphics[width=0.21\textwidth]{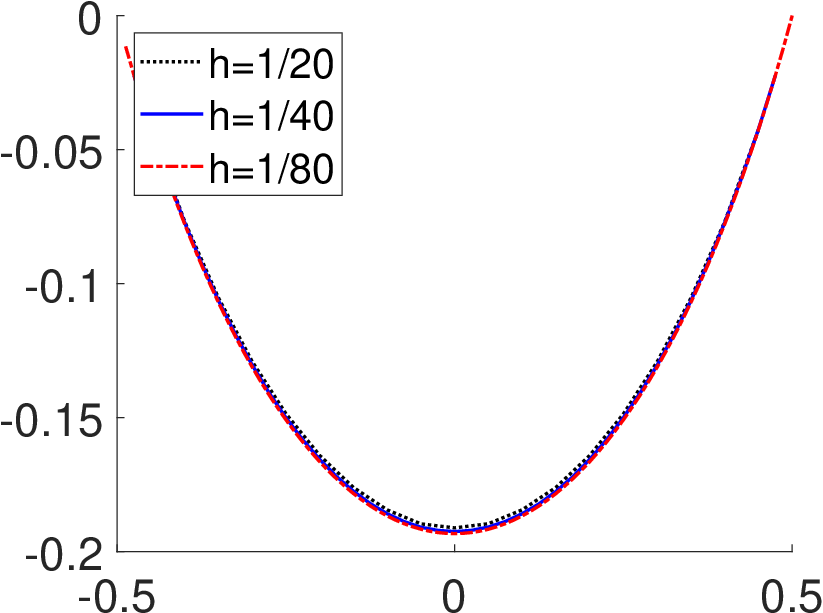}&
		\includegraphics[width=0.21\textwidth]{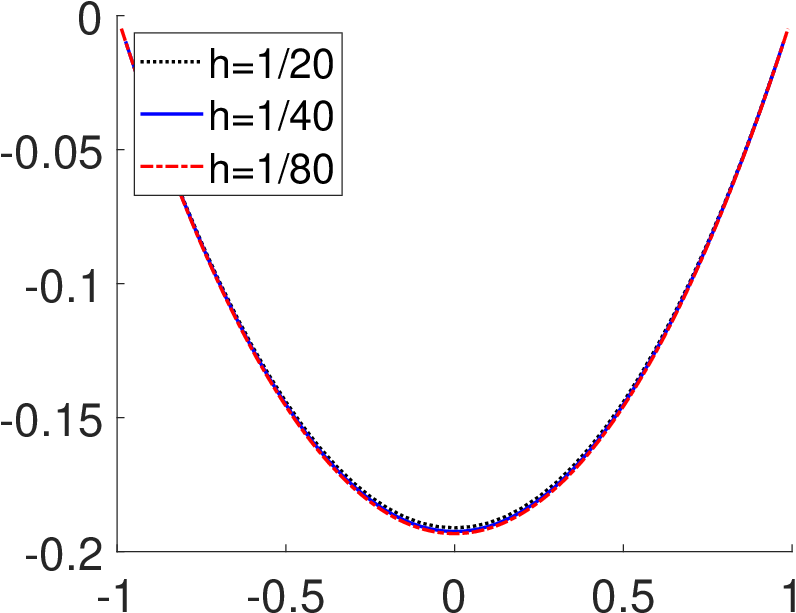}
	\end{tabular}
	\caption{(Test problem (\ref{eq.ex5ellipse}). Scheme (\ref{eq.split1Dis.0})-(\ref{eq.split1Dis.2}).) (a)-(b) Graphs and contour of the numerical solution with $h=1/80$. (c)-(d) Graphs of the restrictions of numerical solutions along $x_1=0$ and $x_2=0$ with $h=1/20$, $1/40$, and  $1/80$. The second-order derivatives are approximated by (\ref{eq.smoothSimp}).}
	\label{fig.ex5ellipse}
\end{figure}

In this example, we solve a problem similar to Example 3 in Section \ref{sec.feng.square}. We consider
\begin{equation}
	K=1/2 \mbox{ in } \Omega,
\end{equation}
and the boundary condition $g=0 \mbox{ on } \partial \Omega$. Different from Example 3, we use the domain being an ellipse:
\begin{equation}
	\begin{cases}
		K=1/2 \mbox{ in } \Omega,\\
		u=0 \mbox{ on } \partial \Omega,\\
		\Omega=\{(x_1,x_2)|x_1^2+4x_2^2\leq 1\}.
	\end{cases}
	\label{eq.ex5ellipse}
\end{equation}
In our experiment, we use scheme (\ref{eq.split1Dis.0})-(\ref{eq.split1Dis.2}) with $\varepsilon=\varepsilon_1=h^2, \Delta t=2h^2$. The second-order derivatives are approximated by (\ref{eq.smoothSimp}). The number of iteration, convergence error and minimum values are shown in Table \ref{tab.ex5}. The graphs and contour of the numerical solution with $h=1/80$ are shown in Figure \ref{fig.ex5ellipse} (a) and (b). The comparison of the restriction of the numerical solution with different $h$ along $x_1=0$ and $x_2=0$ are shown in Figure \ref{fig.ex5ellipse} (c) and (d). Our solution is smooth and convex everywhere.

\section{Conclusion}
\label{sec.conclusion}
In this work, we have proposed two operator splitting/mixed finite-element methods to solve the Dirichlet Minkowski problem in dimension two. Our algorithms are easy to implement since only a system of PDEs is to be solved and the basis functions are chosen to be piecewise linear. When the problem has a classical solution, scheme (\ref{eq.split1Dis.0})-(\ref{eq.split1Dis.2}) using approximation (\ref{eq.double.6})-(\ref{eq.double.7}) for second-order derivatives is first-order accurate while using approximation (\ref{eq.smoothSimp}) it is almost second-order accurate. For an incompatible problem, scheme (\ref{eq.split2Dis.0})-(\ref{eq.split1Dis.2}) can adjust the boundary value of the computed solution to make it compatible with its interior values. Our algorithm can solve the Minkowski problem on arbitrary shaped domains and can also solve problems with singularities in the solution gradient. Our algorithm can be easily extended to high dimensions, which constitutes an ongoing work.

\section*{Acknowledgments}
The preparation of this manuscript has been overshadowed by Prof. Roland Glowinski's passing away. Roland and the authors had intended to write jointly: most of the main ideas were worked out together and the authors have done their best to complete them. In sorrow, the authors dedicate this work to his memory. Roland's creativity, generosity, and friendship will be remembered.

\bibliographystyle{plain}
\bibliography{liu}

\end{document}